\documentclass[11pt]{amsart}

\usepackage{amsthm}
\usepackage{latexsym}
\usepackage[mathscr]{eucal}
\usepackage{ifthen}
\usepackage{mathtools}

\usepackage{amsmath}
\usepackage{amssymb}
\usepackage{graphicx} 
\usepackage{caption} 
\usepackage{subcaption} 

\newtheorem{theorem}{Theorem}[section]

\theoremstyle{definition}

\newtheorem{remark}[theorem]{Remark}

\newcommand{\dt}{\Delta t}

\newcommand{\m}[1]{\mathbf{#1}}

\newcommand{\vy}{\m{y}}
\newcommand{\vb}{\m{b}}
\newcommand{\bt}{\tilde{\vb}}
\newcommand{\bdot}{\dot{\vb}}
\newcommand{\ve}{\m{e}}
\newcommand{\vc}{\m{c}}
\newcommand{\ct}{\tilde{\vc}}

\newcommand{\ep}{\varepsilon}  
\newcommand{\bep}{\vb_\ep}

\newcommand{\mA}{\m{A}}

\newcommand{\Aep}{\mA_\ep}
\newcommand{\mD}{\m{D}}

\newcommand{\mAdot}{\m{\dot{\mA}}}

\numberwithin{equation}{section}

\title[MP-TDRK]{Mixed Precision and Mixed Accuracy Explicit Two-Derivative Runge--Kutta Methods}

\author{Sigal Gottlieb}
\author{Zachary J. Grant}
\author{César Herrera}

\email[Sigal Gottlieb]{sgottlieb@umassd.edu}
\email[Zachary J. Grant]{zgrant@umassd.edu}
\email[César Herrera]{herre125@purdue.edu}

\address[Sigal Gottlieb \& Zachary J. Grant]{Mathematics Department, University of Massachusetts Dartmouth, North Dartmouth, MA 02747.}

\address[César Herrera]{Department of Mathematics, Purdue University, 150 North University Street.
West Lafayette, Indiana 47907-2067}

\thanks{This material is based upon work supported by the National Science Foundation under Grant No. DMS-1929284 while the three authors were in residence at the Institute for Computational and Experimental Research in Mathematics in Providence, RI, during the ``Empowering a Diverse Computational Mathematics Research Community'' program.
SG's research was supported in part by AFOSR Grant No. FA9550-23-1-0037, DOE  Grant No. DE-SC0023164 Subaward RC114586,
and Mass Dartmouth’s Marine and Undersea Technology (MUST) Research Program funded by the ONR Grant No. N00014-20-1-2849.
ZJG's research was supported in part by  DOE grant  No. DE-SC0023164 Subaward RC114586 and AFOSR Grant No. FA9550-23-1-0037.
}

\begin{document}
  \maketitle
  
\begin{abstract}   
Mixed precision Runge--Kutta methods have been recently developed and used
for the time-evolution of partial differential equations.
Two-derivative Runge–Kutta schemes may offer enhanced stability and accuracy properties 
compared to classical one-derivative methods, 
making them attractive in a wide variety of problems. 
However, their computational cost can be significant, 
motivating the use of a mixed-precision paradigm  that employs different floating-point precisions 
for different function evaluations to balance efficiency and accuracy.
To ensure that the perturbations introduced by the low precision computations 
do not destroy the accuracy of the solution, we need to understand how these
perturbation errors propagate.
We extend the numerical analysis  mixed precision framework previously developed 
for Runge--Kutta methods  to characterize the propagation of the perturbation errors 
arising from mixed precision computations in explicit and implicit two-derivative Runge--Kutta methods.
We use this framework for analyzing  the order of the perturbation errors, and 
for designing new methods that are less sensitive to the effect of the low precision
computations.  Numerical experiments on linear and nonlinear representative PDEs, 
demonstrate that appropriately designed mixed-precision two-derivative Runge--Kutta methods 
achieve the predicted accuracy.
  \end{abstract}

\section{Overview}
In numerically solving a partial  differential (PDE) equation of the form
\begin{eqnarray}
\label{pde}
	U_t = f(U, U_x, U_{xx}),
\end{eqnarray}
we often adopt a method of lines (MOL) approach where the spatial derivatives on the 
right hand side are first approximated by some numerical approach 
(e.g. finite differences, finite elements, spectral methods), to obtain the 
semi-discretized system of equations
\begin{eqnarray}
\label{ode}
u_t = F(u),
\end{eqnarray}
(where $u$ is a vector of approximations to $U$). This system can then be evolved forward using a high order
time stepping method, for example a Runge--Kutta method, 
\begin{subequations} \label{RK}
\begin{eqnarray} 
    y^{(i)} & = & y_n + \dt \sum_{j=1}^{s} a_{ij} F(y^{(j)})   \\
    y_{n+1} & = & y_n + \dt \sum_{i=1}^s b_{i} F(y^{(i)}) .
\end{eqnarray}
\end{subequations}
or a  Taylor-series (multiderivative) method,
\begin{eqnarray} \label{TS}
    y_{n+1} & = & y_n   
    +\sum_{i=1}^\ell b_{i} \dt^i \frac{d^i F}{dy^i}(y_n).
\end{eqnarray}
In this work we focus on a combination of these approaches:
multi-stage multi-derivative methods. In particular, we 
explore two-derivative Runge--Kutta methods.

The idea of enhancing  Runge--Kutta methods with additional derivatives was proposed in the 
1950s and 1960s in \cite{Tu50,StSt63}, and  multistage multiderivative time integrators 
for ordinary differential  equations were studied in  \cite{shintani1971,shintani1972,KaWa72,KaWa72-RK, mitsui1982,ono2004, tsai2010}. 
More recently, two-derivative Runge--Kutta  methods were applied to the time evolution of PDEs
 \cite{sealMSMD2014,tsai2014,LiDu2016a,LiDu2016b,LiDu2018}. In this work we consider
two-derivative Runge--Kutta  methods for 
the time-evolution of PDEs of the form
\eqref{pde}, combining this approach with mixed
precision computation meant to improve efficiency while retaining accuracy.

Modern computing architectures have resulted in recent emphasis on mixed precision computations. 
This is because the dominant cost is no longer the computational arithmetic (measured in FLOPS), 
but other factors such as  memory bandwidth, data movement, and energy consumption. Lower precision is often advantageous because it significantly 
reduces these computational costs. However, 
higher precision is often necessary to obtain stable and accurate simulations. 
Mixed precision approaches have become very popular in numerical linear algebra
\cite{perf1,perf5,perf6}. More recently, mixed precision approaches for the numerical solution of 
time-dependent PDEs have been proposed  and shown to be accurate and efficient
\cite{Grant2022,croci2022,BalosRobertsGardner2023,marot2025, MPWENO,Burnett1, Burnett2}. 
These approaches use lower precision when possible, and combine them with higher precision when needed. 
The ideal scenario is a  combination of  the efficiency of low precision
and the accuracy of high precision that allows us to enjoy the best of both worlds. 
In this work we  extend the numerical analysis framework 
for mixed precision Runge--Kutta methods introduced in \cite{Grant2022}
to mixed precision two-derivative Runge--Kutta methods.
This approach allows us to  understand the effect of 
the perturbation errors induced by the low precision computations, and to
design methods that mitigate their impact on the overall solution.

We note that although our primary motivation is mixed precision simulations,
this framework applies to a variety of other scenarios. Whenever a part of the
simulation is computed with less accuracy, the framework we describe in
Section \ref{sec:framework} can be adapted to understand the propagation of the errors.
This framework explains the overall impact of computing only the 
 second order derivative with less accuracy. 
 In fact, in many mixed accuracy methodologies the framework we describe in Section \ref{sec:framework}
 simplifies due to the smoothness of the perturbation introduced, which allows cancellation errors
which cannot be relied upon when the perturbation is non-smooth, such as  in the mixed precision case.

The structure of the paper is as follows: in Section 
\ref{sec:back1} we review the current state of  mixed precision methods for the numerical solution of time-dependent PDEs. In Section \ref{sec:MPRK} we discuss the theoretical framework introduced in \cite{Grant2022}
for mixed precision Runge--Kutta methods.
In Section \ref{sec:MPexp} we briefly discuss the 
mixed precision techniques in exponential time integration methods introduced in \cite{BalosRobertsGardner2023},  
and in Section \ref{sec:MPRKC} the 
mixed precision Runge-Kutta--Chebyshev methods of \cite{croci2022}.
In Subsection \ref{sec:MPWENO} we review the mixed precision WENO finite difference method introduced in
\cite{MPWENO}.
In Section \ref{sec:back2} we present two-derivative methods and their recent use in hyperbolic PDEs. 

In Section \ref{sec:framework} we present the mixed precision two-derivative Runge--Kutta framework, 
following the approach in \cite{Grant2022}. We use a simplified framework that assumes that only
computations of the second derivative are performed in low precision, and use 
this framework to establish  order conditions for the low precision perturbation.
In Section \ref{sec:MPTDRKmethods} we present third, fourth, fifth, and sixth order mixed 
precision methods with various perturbation orders.
In Section \ref{sec:NumResults} we present the mixed precision numerical simulations that show that the numerical results are as predicted by the theory. Finally, in Section \ref{sec:conclusions} we present our conclusions 
and discuss possible future directions of this approach.

\section{Background: mixed precision methods for  PDEs} \label{sec:back1}
The use of mixed precision for the solution of PDEs is relatively new.
In this section we review the existing work on mixed precision 
methods used in the context of PDEs.

\subsection{Mixed Precision Runge--Kutta methods} \label{sec:MPRK}
A mixed accuracy framework for Runge--Kutta methods was introduced 
in \cite{Grant2022} with an emphasis on  designing mixed-precision 
implicit Runge–Kutta methods. This work utilizes lower-precision arithmetic 
in the implicit computation of a multi-stage time integrator without 
destroying  the overall accuracy of the solution.  The mechanism 
behind this approach takes advantage of the multi-stage nature of 
the scheme to damp out the errors from the low-precision computation.
The key idea is that the implicit solves (or the expensive part of 
the implicit solve)  can be evaluated in low precision such as  half 
(\texttt{Float16})
or single (\texttt{Float32}) precision,
while the explicit computations are computed 
in high precision, leading to a highly accurate final update.

As a simple example, consider the   diagonally  implicit Runge--Kutta method 
\begin{subequations} \label{DIRK}
\begin{eqnarray} 
    y^{(i)} & = & y_n + \dt \sum_{j=1}^{i} a_{ij} F(y^{(j)})   \\
    y_{n+1} & = & y_n + \dt \sum_{i=1}^s b_{i} F(y^{(i)}) ,
\end{eqnarray}
\end{subequations}
where each stage is computed implicitly, and the final reconstruction is explicit.
To make the implicit stages $ y^{(i)}$ cheaper to invert, we  replace $F(y^{(i)})$ with a lower precision function $F_{\varepsilon}(y^{(i)})$. This is harmonious with the numerical linear algebra approaches that
iteratively solve the nonlinear system using a mixture of  precisions \cite{CTKelley2022}.
The mixed precision method then becomes:
\begin{subequations} \label{MPDIRK}
\begin{eqnarray} 
    y^{(i)} & = & y_n + 
    \dt \left( \sum_{j=1}^{i-1} a_{ij} F(y^{(j)})
    + a_{ii} F_{\varepsilon}(y^{(i)}) \right) \\
    y_{n+1} & = & y_n + \dt \sum_{i=1}^s b_{i} F(y^{(i)}) .
\end{eqnarray}
\end{subequations}
This strategy introduces a perturbation  into the internal stages. 
The key to mixed precision
methods is that this perturbation must be controlled 
by computing (and when needed re-computing) the function evaluations
in high precision  so that the overall solution is minimally impacted by it.

Grant \cite{Grant2022} introduced a general theoretical additive perturbation 
framework showing how perturbations from low precision can be controlled. 
Within this framework, specific order conditions are derived that track the propagation 
of the perturbation errors introduced by the lower precision computations, 
and enable the characterization of which operations   may be safely computed in lower precision.

The mixed precision implementation of the fully implicit Runge--Kutta method \eqref{RK}
can be re-written as 
\begin{subequations} 
\begin{align}
\label{Perturbed Methoda}
\vy &= u^n \ve  + \dt
\tilde{\mA}_{ij} F(\vy) +  \ep \dt \Aep \tau(\vy)\\
u^{n+1}&= u^n + \dt \tilde{\vb}_{j} F(\vy) + \ep  \dt  \bep \tau(\vy) \label{Perturbed Methodb}
\end{align}
\end{subequations}
where $\tilde{A}= \mA+ \mA_\ep$ and $\tilde{\vb} = \vb+ \bep$,
and 
\[  \tau(y) = \frac{1}{\varepsilon} \left( F(y) - F_{\varepsilon}(y)\right). \]

Analyzing the scheme in this form allows us to use an additive 
B-series representation to track the evolution of $ u$ 
as well as its interaction with  the perturbation  function $\tau$. 
For example, a second order expansion is:
\begin{align*} 
u^{n+1} &= u^n +  \dt \bt F(\vy) + \ep \dt \bep \tau(\vy) 
\nonumber \\
&= u^n +  \dt \bt F\left(u^n \ve + \dt \tilde{\mA} F(\vy)
+ \ep  \dt  \Aep \tau(\vy) \right) + \ep \dt \bep \tau(\vy) 
\nonumber \\
&= u^n +  \dt \bt F\left(u^n \ve \right) + 
\dt^2 \bt \tilde{\mA}  F_y\left(u^n \right) F(\vy)
+ \ep  \dt^2  \bt \Aep  F_{y} \left(u^n\right) \tau(\vy) 
\nonumber \\
& + \ep \dt \bep \tau(\vy)  + O(\dt^3) + O(\ep \dt^3)
\nonumber 
\end{align*} 
so that
\begin{align} \label{expansion}
u^{n+1} & = \underbrace{ u^n +  \dt \bt \ve F(u^n) + \dt^2\bt \ct F_y(u^n)F(u^n)}_{scheme}  \\
& +   \underbrace{\ep \left( \dt \bep  +  
 \dt^2 F_y(u^n) \bt \Aep  \right)\tau(\vy) }_{perturbation}
\nonumber   + O(\dt^3) + O(\varepsilon \dt^3). \nonumber 
\end{align} 
This expansion shows two sources of error: those of the scheme and those of the perturbation.

In the case of a smooth perturbation, we can rely on cancellations, and even expand the $\tau$ term,
and obtain more relaxed conditions. However, in the case where the perturbation is non-smooth, 
such as mixed precision,
we  cannot rely  on cancellations, or expand the $\tau$ terms. In this case, if we wish to eliminate the 
$\ep \dt$ term, it is not enough to require that
$\bep \ve =0$, rather we require that each element of $\bep$ is zero.


We can easily extend \eqref{expansion} to higher order using the 
additive B-series analysis of the 
scheme. 
The perturbation errors feature a combination of 
the scheme coefficients and the perturbation coefficients. Expanding the terms
to third order in $\dt$ and first order in $\epsilon$ we obtain:
\[ \begin{array}{lll} \hline 
\mbox{Terms involving } \; \; \; \;  & \mbox{Terms involving } & \mbox{scheme} \\ 
\mbox{ $\ep$ and $\dt$} \; \; \; \;  & \mbox{ $F$ and $\tau$} & \mbox{coefficients} \\ \hline
& & \\
\ep \dt \; \; \; \; &  \tau(u^n)      & \bep \ve \\ 
\ep \dt^2 & F_y(u^n)\tau(\cdot)    & \bt \Aep \ve\\ 
\ep \dt^3 & F_y(u^n)F_y(u^n)\tau(\cdot) &  \bt \tilde{\mA} \Aep  \ve \\
\ep \dt^3  & F_{yy}(u^n)(F(u^n),\tau(\cdot))    
&  \bt(\ct \cdot \Aep \ve)  \\
 \end{array} \]
From these terms, the  consistency conditions for the scheme and the perturbation conditions can be easily defined, even if the perturbation $\tau$ is poorly behaved, as is expected in mixed precision. 
For a scheme to achieve order $O( \ep \dt^{m+1})$ we require :\\
\begin{tabular}{ll}
$m\geq 1$ & $|\bep| \ve=0$ \\
$m\geq 2$  & 
  $|\bt| | \mA_\ep| \ve = 0 $\\
$m\geq 3$  & 
 $| \bt| | \Aep| |\tilde{\mA}| \ve  = 0$ , \; \;
 $|\bt| |\tilde{\mA}| | \mA_\ep| \ve =0$. \\ 
\end{tabular} \\
We observe that we must have each element of $\bep $ be zero,
and the other conditions depend on which stages are excluded from the final row reconstruction.

Grant's order conditions in \cite{Grant2022} 
show how to design a method
to ensure that the low-precision perturbation does not ruin the accuracy of the method. In the simplest and most intuitive
sense, we see that a diagonally implicit  method of the form \eqref{MPDIRK} above introduces the perturbation into the implicit stage only, 
so that each stage has a perturbation error of
$\dt \varepsilon$, and this is further mitigated at the final stage by another $\dt$, so that at each step we have a perturbation error of
$\dt^2 \varepsilon$. Growing to a final time of $T_f$, we expect $O(\frac{1}{\dt})$ steps and so a final perturbation error of $\dt \varepsilon$.
Explicit \cite{Grant2022, Burnett2} or inexpensive implicit \cite{Grant2026} corrections can then be used further to  improve the perturbation error. 

Performance analysis of methods developed or analyzed by  the framework 
introduced in \cite{Grant2022} was performed in \cite{Burnett1,Burnett2} for small problems, and in \cite{AlSayedAli2025} on larger ODEs where the implicit problem is solved in low and mixed precision. In 
\cite{Dravins2025} the methods of \cite{Grant2022} were tested 
with state-of-the-art software using the Ginkgo library.
\begin{remark} {\em
    We note that while mixed precision has only recently been used in this context, the implicit stages are rarely computed exactly and some perturbation  is introduced  whenever the implicit stage is approximated by  some iterative procedure such as  Newton's iteration.
For this reason, perturbation analysis is always important for implicit methods; however, the perturbations introduced by mixed precision are non-smooth, and so more difficult to handle.
 The mixed precision order conditions introduced in \cite{Grant2022} allow us to deal with any type of perturbation, whether smooth (such as linearization) or non-smooth (such as rounding error in low precision). }
\end{remark}
\subsection{Mixed precision exponential time stepping methods} \label{sec:MPexp}
Mixed precision approaches have been most heavily developed in the context of  numerical linear algebra methods. It is natural, then, that numerical methods for time-evolution that have the closest relationship to numerical
linear algebra will take advantage of the mixed precision paradigm.
In \cite{BalosRobertsGardner2023}, the mixed precision paradigm was
leveraged for efficiency in solving  an advection-diffusion-reaction PDE
using  exponential time integration methods. 
Two approaches were studied in conjunction with  exponential time integrators. 
First, the authors reformulated the exponential Rosenbrock-Euler method in a way that allows for low precision computations of the matrix exponential which improves efficiency while improving accuracy over the low precision version. Next, the authors use a mixed precision 
Krylov approach  for computing matrix exponentials, which was more efficient than the full precision approach and more accurate than the low precision approach.

\subsection{Mixed precision explicit Runge--Kutta--Chebyshev methods}
\label{sec:MPRKC}
 While the methods developed  \cite{Grant2022} and the follow-up work in \cite{Burnett1,Burnett2,Dravins2025,AlSayedAli2025} 
 focused on implicit Runge--Kutta methods, where the expensive 
 implicit solves dominate runtime, 
 this is not the only avenue in which  mixed precision has significant potential to
 accelerate computations without significantly degrading accuracy. 

When dealing with ``stiff'' differential equations,
where the eigenvalues are on the negative real axis, 
the solution decays rapidly, 
necessitating a careful choice of the numerical method. 
Typically, explicit methods are inexpensive but
require very small time-steps when solving the such equations. 
However, Runge--Kutta--Chebyshev methods are explicit first and second order 
Runge--Kutta type methods with a large number of stages $s$, 
that feature  very large  linear stability regions along the 
negative real axis. 
Just like the Gauss-Chebyshev-Lobatto quadrature methods, these
Runge--Kutta methods use the recurrence relation definition 
of the Chebyshev polynomial to create a low storage formula that results 
in a second order approximation to the solution of an ODE system.
If a large linear negative-real-axis stability region is desired, this 
method may have a very large number of stages $s$. Furthermore,
this method is built recursively at each time-step. These two factors make these
explicit Runge--Kutta methods computationally costly.

However, since the polynomial is built recursively, one can achieve high order accuracy 
by computing most of the stages  in low precision, with only a few stages  computed 
in high precision.
Using this insight, Croci \& de Souza \cite{croci2022} suggested   
a mixed precision approach that exploits a clever re-writing 
of these methods to allow for an inexpensive mixed precision implementation
to design explicit Runge--Kutta Chebyshev methods with a mixed precision paradigm.
This work showed that while careless use of low precision may break convergence,
carefully chosen high-precision function 
evaluations allow for  order-preserving mixed precision that is computationally efficient.

\subsection{Mixed Precision WENO} \label{sec:MPWENO}
 Another notable use of mixed precision in time-dependent simulation of a hyperbolic problem was in \cite{MPWENO}, where  a novel mixed-precision weighted essentially non-oscillatory (WENO) method was developed and used for solving the Teukolsky equation.
 This hyperbolic PDE arises when modeling perturbations of Kerr black holes. 
 In these simulations, we require very long-time evolution, in which rounding error may build up. 
 For this reason, high order discretizations  and quad precision (\texttt{Float128})
 are needed to accurately evolve the solution. In \cite{MPWENO}, the authors showed that 
 while quad precision simulation is needed for the long time evolution, not all parts of the WENO algorithm need be performed in quad precision. In particular, the expensive nonlinear 
 computation of the WENO weights was performed in reduced  precision, 
while requiring that the WENO weights sum to one in quad precision.
  The optimized and accelerated WENO code took seven days to run in quad precision. 
The mixed precision implementation  took only two days to run, 
a significant speedup factor of 3.3x, while producing equally accurate results.

\section{Background: two derivative Runge--Kutta 
methods} \label{sec:back2}

Two derivative Runge--Kutta (TDRK) methods  were applied to 
the time evolution of PDEs in  \cite{sealMSMD2014,tsai2014,LiDu2016a,LiDu2016b,LiDu2018}.
When applied to a system of ODEs \eqref{ode} resulting from a semi-discretization of a PDE \eqref{pde}, 
a TDRK  can be written as:
\begin{eqnarray}
\label{TDRK}
y^{(i)} & = &  u^n +  \dt \sum_{j=1}^{s}  a_{ij} F(y^{(j)}) +
\dt^2   \sum_{j=1}^{s} \dot{a}_{ij} \dot{F}(y^{(j)}) , \; \; \; \; i=1, . . ., s \\
 u^{n+1} & = &  u^n +  \dt \sum_{j=1}^{s} b_{j} F(y^{(j)}) +
\dt^2  \sum_{j=1}^{s} \dot{b}_{j} \dot{F}(y^{(j)})  , \nonumber
\end{eqnarray}
or, in matrix form
\begin{eqnarray}
\label{TDRKmatrix}
\vy & = &  u^n +  \dt \mA F(\vy) +
\dt^2 \mAdot \dot{F}(\vy)  \\
 u^{n+1} & = &  u^n +  \dt \vb F(\vy) + \dt^2 \bdot \dot{F}(\vy)  \nonumber
\end{eqnarray}
where
$ \mA_{ij}  =  a_{ij} $ , 
$  \mAdot = \dot{a}_{ij}$, 
 and $\vb$ and $\dot{\vb}$ are vectors with the elements $b_{j}$ and $\dot{b}_{j}$, respectively.
The order conditions for this method up to sixth order
are in given in Appendix A of \cite{TSpaper}.
In this work we wish to focus only on explicit methods, 
so that the sums at the $i$th stage of  \eqref{TDRK}  go up to $i-1$ rather than $s$,
and correspondingly the matrices $\mA$ and $\mAdot$ are strictly lower triangular.

We note that in Tsai's work \cite{tsai2010} 
a simplified (and less general) version of these methods is given in the form
\begin{eqnarray} \label{Tsai}
y^{(1)} & = &  u^n \\
y^{(i)} & = &  u^n +  \dt a_{i1} F(y^{(1)}) +
\dt^2   \sum_{j=1}^{i-1} \dot{a}_{ij} \dot{F}(y^{(j)}) , \; \; \; i=2, . . ., s  \nonumber \\
 u^{n+1} & = &  u^n +  \dt  F(y^{(1)}) + \dt^2  \sum_{j=1}^{s} \dot{b}_{j} \dot{F}(y^{(j)}). 
\end{eqnarray}
An advantage of this form is the simplified order conditions. 
In some cases this form  may also be advantageous for 
mixed precision purposes.

TDRKs are especially appealing when dealing with  hyperbolic  conservation laws  of the form
\begin{eqnarray}
\label{hyp}
	U_t + f(U)_x = 0,
\end{eqnarray}
where the computation of  $f'(U)$ is generally needed for flux splitting in the spatial discretization. The second derivative computation $\dot{F}$ can re-use the computation of this Jacobian, which makes the TDRK methods natural in this context.

Consider the spatial discretization of \eqref{hyp}, 
where the method-of-lines approach leads to the system 
\[ u_t  = - \mD_x f(u) = F(u) ,\]
for some differentiation matrix  $\mD_x$.
The computation of the second derivative term $\dot{F}$ follows directly from the 
definition  of $F$,
\[ \dot{F} = F(u)_t = F_u u_t = F_u F.\]
In practice  \cite{SDpaper} 
we can use a Lax-Wendroff type approach to compute $U_{tt} \approx \dot{F}$:
\[  U_{tt}  =  - f(U)_{xt}  =  - \left(- f'(U) f(U)_x \right)_x ,\] 
so that 
\[  \dot{F}(u) \approx  - \mD_x \left(f'(u) F(u) \right) ,\] 
where  $f'(u)$ is pre-computed for the flux splitting.
In the next section we develop a mixed precision framework for TDRK methods.

\section{A mixed precision  framework for two-derivative Runge--Kutta methods} \label{sec:framework}

We  consider a mixed precision approach to the 
two-derivative Runge--Kutta method \eqref{TDRK}.
In principal, we can consider either $F$ or $\dot{F}$ (or both) to be in low precision, 
but we currently focus on the case where the computation of the second derivative is performed in low precision, due to the fact that this term enters with a $\dt^2$ multiplication which
naturally reduces the impact of the low-precision perturbation. A more general approach is 
possible, and may be explored in future work.

The low-precision computation of the second derivative is particularly
advantageous where this term is costly to compute and we wish to compute it in a cheaper, albeit less accurate, way. We focus on the mixed precision form 
\begin{eqnarray}
\label{MPTDRK}
y^{(i)} & = &  u^n +  \dt \sum_{j=1}^{s}  a_{ij} F(y^{(j)}) +
\dt^2   \sum_{j=1}^{s} \dot{a}_{ij} \dot{F}_\ep(y^{(j)}) , \; \; \; \; i=1, . . ., s 
\nonumber \\
\; \; \;  \; \; \; u^{n+1} & = &  u^n +  \dt \sum_{j=1}^{s} b_{j} F(y^{(j)}) +
\dt^2  \sum_{j=1}^{s} \dot{b}_{j} \dot{F}_\ep(y^{(j)})  , 
\end{eqnarray}
where $\dot{F}_\ep$ denotes the low precision computation of the second derivative $\dot{F}$. 
In matrix form, this becomes
\begin{eqnarray}
\label{MPTDRKmatrix}
\vy & = &  u^n +  \dt \mA F(\vy) +
\dt^2 \mAdot_\ep \dot{F}_\ep(\vy) \\
 u^{n+1} & = &  u^n +  \dt \vb F(\vy) + 
 \dt^2 \bdot_\ep \dot{F}_\ep(\vy) , \nonumber
\end{eqnarray}
Where the coefficients $\mAdot$ and $\bdot$ in \eqref{TDRKmatrix} are replaced by
$\mAdot_\ep$ and $\bdot_\ep$ to emphasize that they correspond to low precision computations.

We extend the mixed precision framework in \cite{Grant2022} to \eqref{MPTDRKmatrix}, 
first by re-writing the method  in the perturbed form. To distinguish between the behavior of the derivative which should contribute to the accuracy, and the behavior of the perturbation, which we wish to eliminate when possible,
we use the notation 
\[\dot{F}_\ep  = \dot{F} + \varepsilon \dot{\tau}\]
where $\dot{\tau}$ is \emph{not} the derivative of the perturbation $\tau$, but a perturbation of the derivative $\dot{F}$.
Plugging this into \eqref{MPTDRKmatrix} we obtain
\begin{eqnarray}
\label{MPTDRK2}
\vy & = &   u^n +  \dt \mA F(\vy) +
\dt^2 \mAdot \dot{F}(\vy)  +
\ep \dt^2 \mAdot_\ep \dot{\tau}\\
 u^{n+1} & = &  u^n +  \dt \vb F(\vy) + 
 \dt^2 \bdot \dot{F}
 + \ep \dt^2 \bdot_\ep \dot{\tau} ,\nonumber
\end{eqnarray}
where for clarity we use the coefficients $\mAdot$ and $\bdot$ attached to $\dot{F}$
and the coefficients $\mAdot_\ep$  $\bdot_\ep$ attached to $\dot{\tau}$, even though in our case
$\mAdot = \mAdot_\ep$ and $\bdot = \bdot_\ep$. 

We can Taylor expand this form to 
obtain the following terms, and match them to 
the desired order:\\
\begin{center}
\begin{tabular}{llll} 
$O(\dt)$: & $\vb \ve = 1$ &  \\
$O(\dt^2)$: &  
$ \vb \vc +  \dot{\vb} \ve = \frac{1}{2}$ \; \; \; \; & 
$O(\ep \dt^2)$:  
& $ \left| \dot{\vb}_\ep \right| \ve = 0 $ \\
$O(\dt^3)$: &  $ \vb \mA \vc  + \vb \mAdot \ve + \dot{\vb} \vc= \frac{1}{6} $ \; \; \; \; &
$O(\ep \dt^3)$: & 
$ \left| \vb \right| \left| \mAdot_\ep \right| \ve = 0$ \\
$O(\dt^3)$: & $ \vb  \vc^2 + 2 \dot{\vb} \vc = \frac{1}{3} $ &
$O(\ep \dt^3)$: & $ \left| \dot{\vb}_\ep \right| \left|\vc\right| = 0$ \\ 
\end{tabular}
\end{center}
where $\vc^2$ is understood as element-wise squaring, and the absolute value $| \cdot |$ 
is understood element-wise both for matrices and vectors.
In this work we limit the perturbation conditions to third order.
Expanding to higher order terms will produce more order and perturbation conditions.


\smallskip

\noindent{\bf First order perturbation error:}
All methods that satisfy the formulation \eqref{MPTDRKmatrix} where all the evaluations of $F$ 
are in full precision, have the property that   the perturbation enters  as $\ep \dt^2$ at each time-step. 
Once again, evolving to a final time requires $O(\frac{1}{\dt})$ steps, 
so that the perturbation error builds up to $E_{pert} = O(\ep \dt)$. 
Thus, all methods of the form \eqref{MPTDRKmatrix} have overall final-time
error
\[ E = E_{method} + E_{pert} = O(\dt^p) + O(\ep \dt).\]
We expect to see an order of $p$ as long if $\ep << \dt^{p-1}$, and first order as 
$\dt$ gets smaller.

\smallskip

\noindent{\bf Second order perturbation error:}
To improve the errors resulting from the low precision perturbation,
we can design a TDRK method so that the perturbation condition 
corresponding to $\ep \dt^2$ is 
satisfied:
\begin{eqnarray}\label{pert2}
     \left| \dot{\vb}_\ep \right| \ve = 0, \; \; \rightarrow  \left(\dot{\vb}_\ep\right)_j = 0, \; \; \forall j.
\end{eqnarray}
This condition means that the second derivative terms computed in low precision do not enter 
in the reconstruction stage.
A method of the form \eqref{MPTDRKmatrix}  that additionally 
satisfies this perturbation condition, will have final time-error
\[ E = E_{method} + E_{pert} = O(\dt^p) 
+ O(\ep \dt^2).\]
We expect to see  order $p$ convergence   if 
$\ep << \dt^{p-2}$, and second order
for the regime where $\dt$ gets smaller.

\smallskip

\noindent{\bf Third order perturbation error:}
Continuing with this approach, we can zero out the perturbation errors 
that correspond the $\ep \dt^3$, by additionally requiring
$ \left| \dot{\vb}_\ep \right| |\mA | \ve = 0$
and 
$ |\vb| \left| \mAdot_\ep \right| \ve = 0$. The first of these conditions
is automatically satisfied by $\left(\bdot_\ep\right)_j = 0, \; \; \forall j$, so we 
need only impose 
\begin{eqnarray} \label{pert3}
\; \; \; \sum_{i=1}^s |b_i| \sum_{j=1}^s \left| \left(\mAdot_\ep\right)_{ij} \right|  = 0 .
\end{eqnarray} 
This condition ensures that the method does not include in the final reconstruction any stage that has a perturbed value,
which  gives us  a perturbation error of 
$O(\ep \dt^4)$. At the final time
we expect an error of
\[ E = E_{method} + E_{pert} = O(\dt^p) 
+ O(\ep \dt^3).\]
We expect to see an order of $p$ as long if 
$\ep << \dt^{p-3}$, and third order as $\dt$ gets smaller.

\begin{remark} \label{rmk2}
{\em It is important to note that the value of $\ep$ is not the machine zero associated with the low precision.
    While  $\ep$ is clearly connected with the number
    of digits retained in the lower precision, it
    also depends on the number of computations to evaluate the operator, 
    and how they are performed. 
    For example, if we square a low precision matrix of dimensions $N_x \times N_x$, we expect that the
    $N_x^2$ multiplications and $N_x$ additions performed to evaluate each low precision computation would result in a much larger $\ep$ than one would naively expect from the low precision machine zero! On the other hand, performing 
    this matrix multiplication in high precision and then casting
    it to lower precision for use at each stage will allow
    for a smaller $\ep$, but would increase the computational cost.
    We see in the numerical examples below (Figure \ref{fig:linadv3p}) the impact of the computation of the low precision operator on the value of $\ep$. } 
\end{remark}

\section{Explicit mixed precision TDRKs} \label{sec:MPTDRKmethods}

In this section  focus on explicit mixed precision TDRKs. 
We  provide explicit mixed precision methods that have final 
time errors
\[ E = O(\dt^p) + O(\ep \dt^m) \]
of orders $p=3,4,5,6$,  and various perturbation orders $m$. 
In designing the methods we prioritize simple
rational coefficients, which are beneficial for implementation in quad precision.
We present the linear stability regions of these methods, as well
as the impact of a perturbation on these stability regions.
In Section \ref{sec:NumResults} we will verify that these methods
obtain the design order of accuracy and perturbation.

\subsection{Third order methods} \label{sec:MPTDRK3p}

We begin with a method  found in \cite{tsai2010}:

\smallskip

\noindent{\bf TDRK method TDRK2s3p1e:}
The  two-stage third order method  
\begin{align}  \label{MP2s3p1e}
  y^{(1)}    &= u^n + \Delta t F(u^n) + 
  \frac{\Delta t^2}{2} \dot{F}_\ep(u^n) \nonumber \\
    u^{n+1} &= u^n + \Delta t F(u^n) + 
    \frac{\Delta t^2}{6}( 2 \dot{F}_\ep(u^n) + \dot{F}_\ep(y^{(1)}))
\end{align}
will have final time error 
\[ E = O(\dt^3) + O(\ep \dt).\]

This method takes  Tsai's form \eqref{Tsai} so that  we only need to compute the 
high precision value $F(u^n)$ once per time-step, which also provides computational savings. 
We can improve the perturbation error by requiring that the final reconstruction does 
not include any second derivative terms, but then Tsai's form is not possible.

\begin{figure}[ht] 
    \centering
    \begin{subfigure}[c]{0.325\textwidth}
        \centering
        \includegraphics[width=\textwidth]{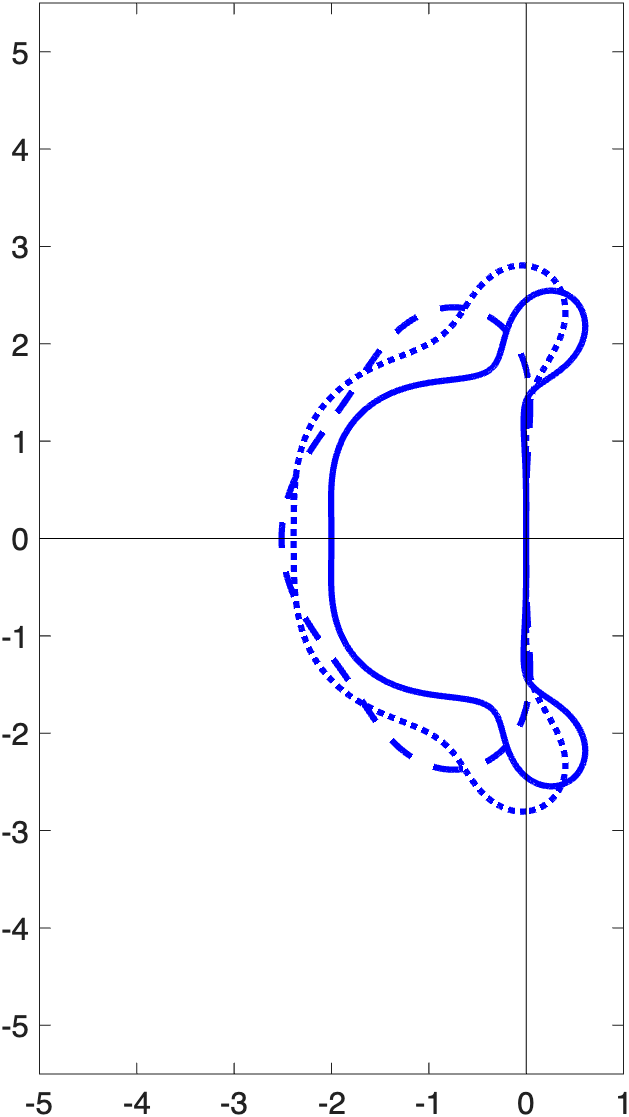}
    \end{subfigure}
    \hfill
      \begin{subfigure}[c]{0.325\textwidth}
        \centering
\includegraphics[width=\textwidth]{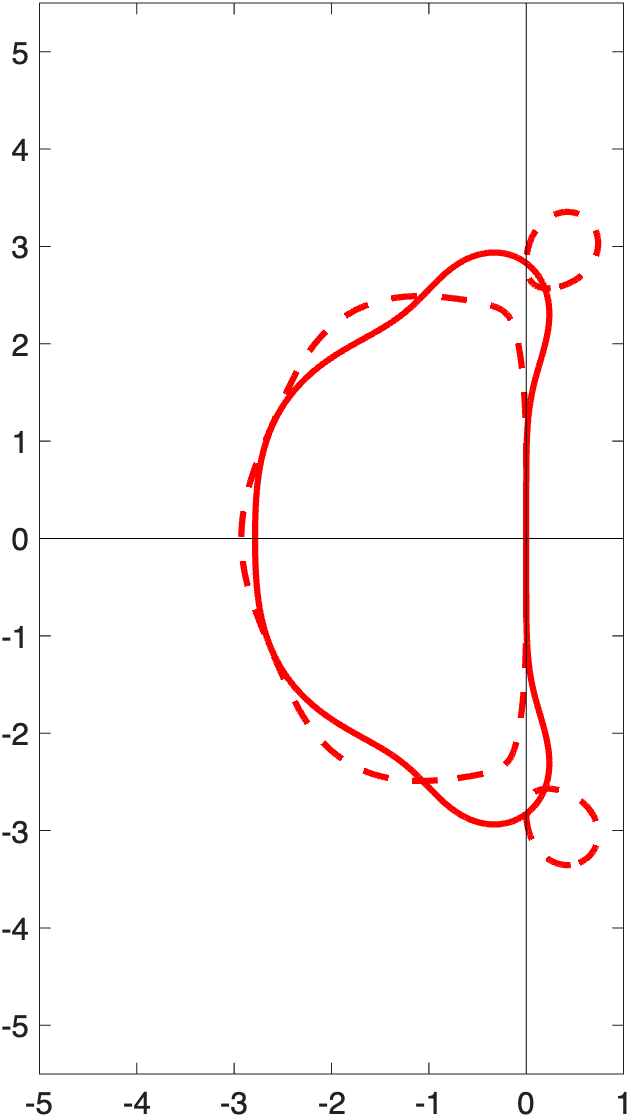} 
    \end{subfigure}
    \hfill
    \begin{subfigure}[c]{0.325\textwidth}
        \centering
\includegraphics[width=\textwidth]{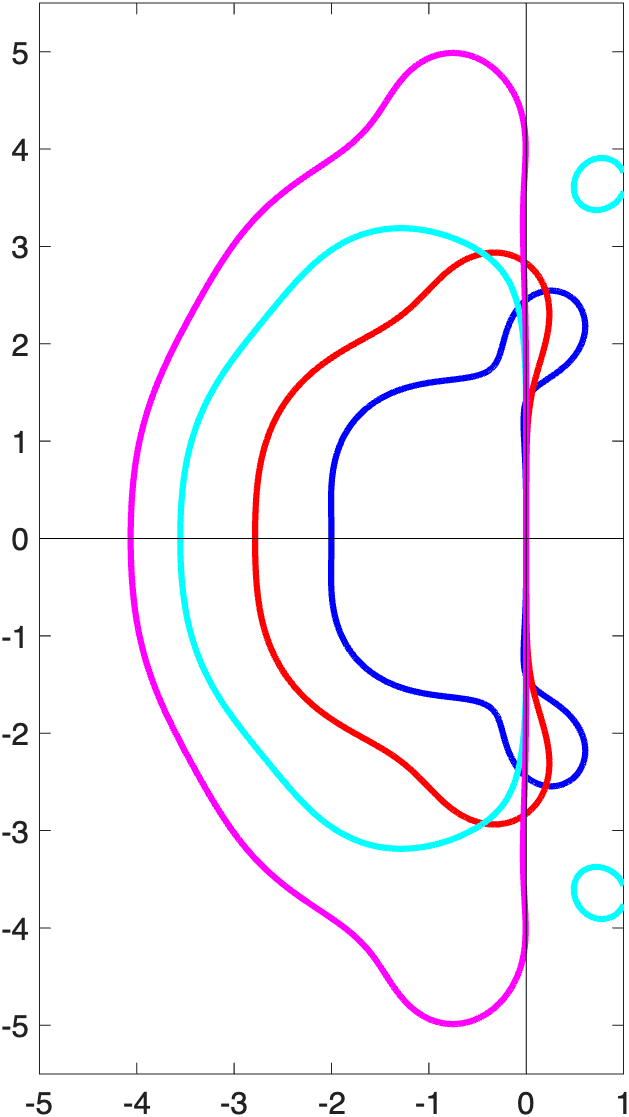} 
    \end{subfigure}
    \hfill
     \caption{Linear stability regions for the MP-TDRK methods with error
     $O(\dt^p)$ and perturbation error $O(\ep \dt^m)$.
     Third order in blue, fourth order in red, fifth order in cyan, 
     sixth order in magenta. Solid lines are perturbation error $m=1$,
     dashed lines are $m=2$, and dotted line $m=3$.
     Left: Third order methods. Center: Fourth order methods.
     Right: All methods with     perturbation order $m=1$.}
     \label{fig:stability}
\end{figure}

\smallskip

\noindent{\bf TDRK method TDRK2s3p2e:}
This two-stage third order method satisfies Condition \eqref{pert2}: 
\begin{align}  \label{MP2s3p2e}
y^{(1)}    &= u^n + \frac{2}{3} \Delta t F(u^n) 
    + \frac{2}{9} \Delta t^2 \dot{F}_\ep(u^n) \nonumber \\
u^{n+1} &= u^n + \frac{1}{4} \Delta t F(u^n) 
+ \frac{3}{4} \Delta t F(y^{(1)}) 
\end{align}
This method will have a final time error of 
\[ E = O(\dt^3) + O(\ep \dt^2).\]

\begin{figure}[t] 
    \centering
    \begin{subfigure}[c]{0.27\textwidth}
        \centering
\includegraphics[width=\textwidth]{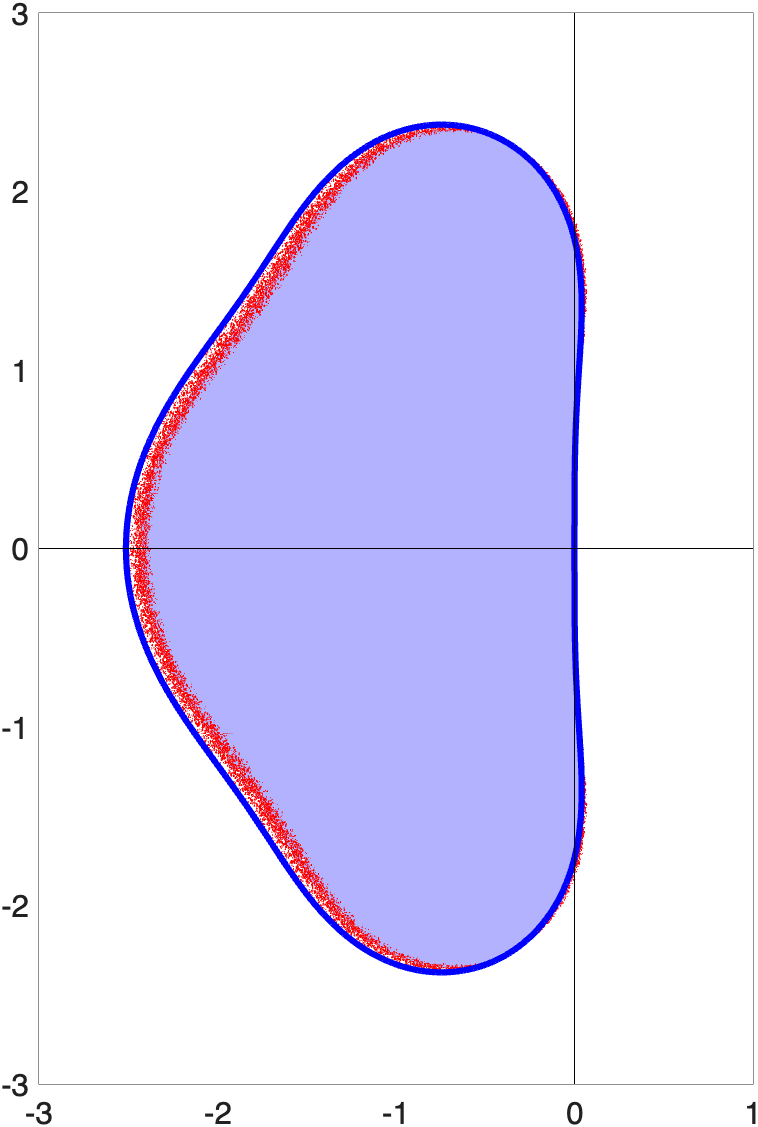}
\includegraphics[width=\textwidth]{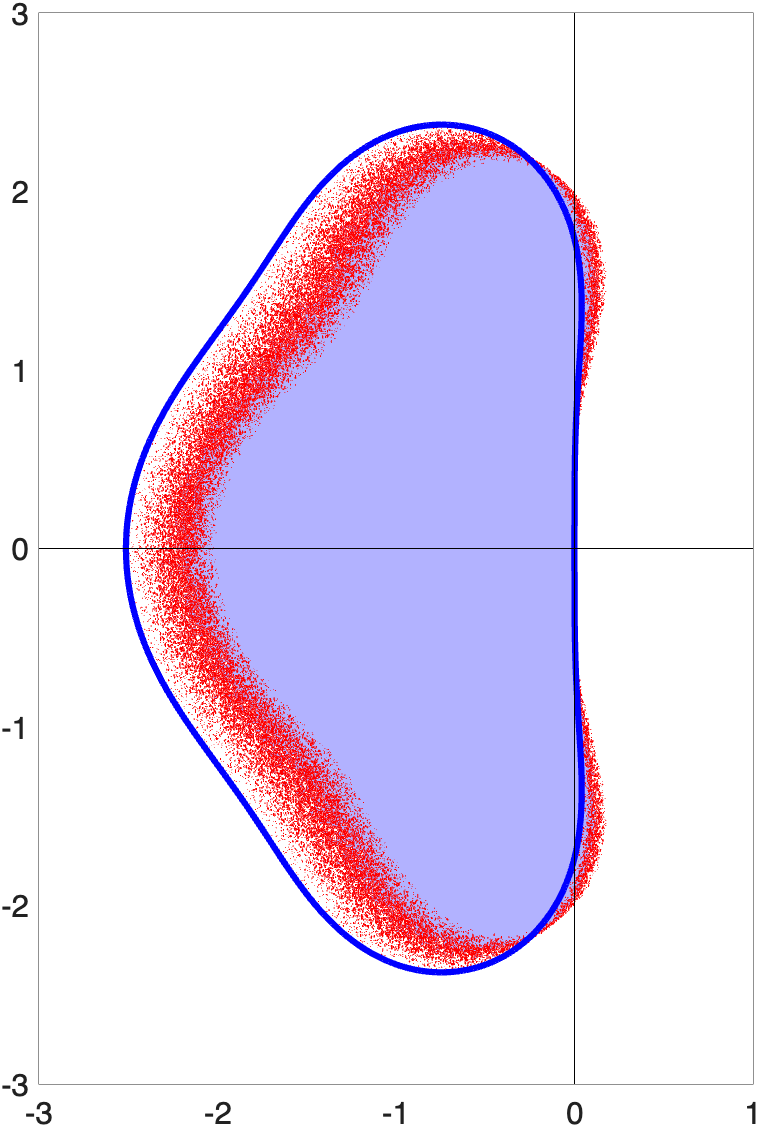}
    \end{subfigure}
      \begin{subfigure}[c]{0.27\textwidth}
        \centering
 \includegraphics[width=\textwidth]{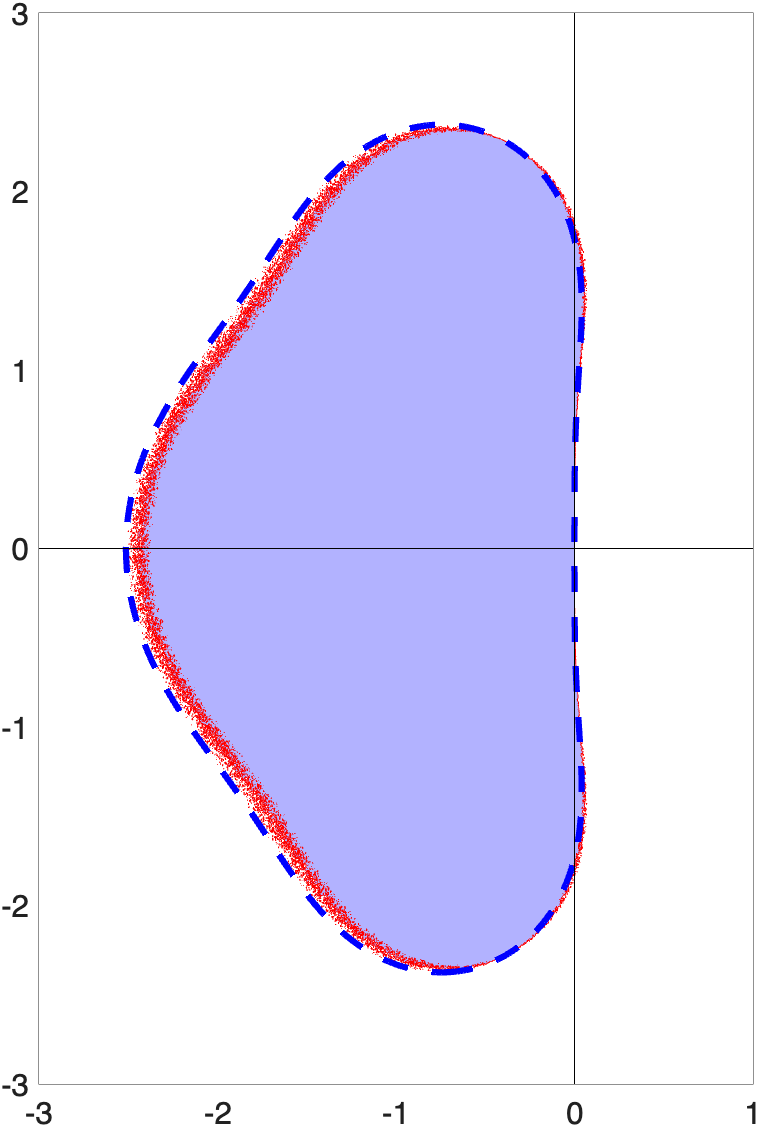}
  \includegraphics[width=\textwidth]{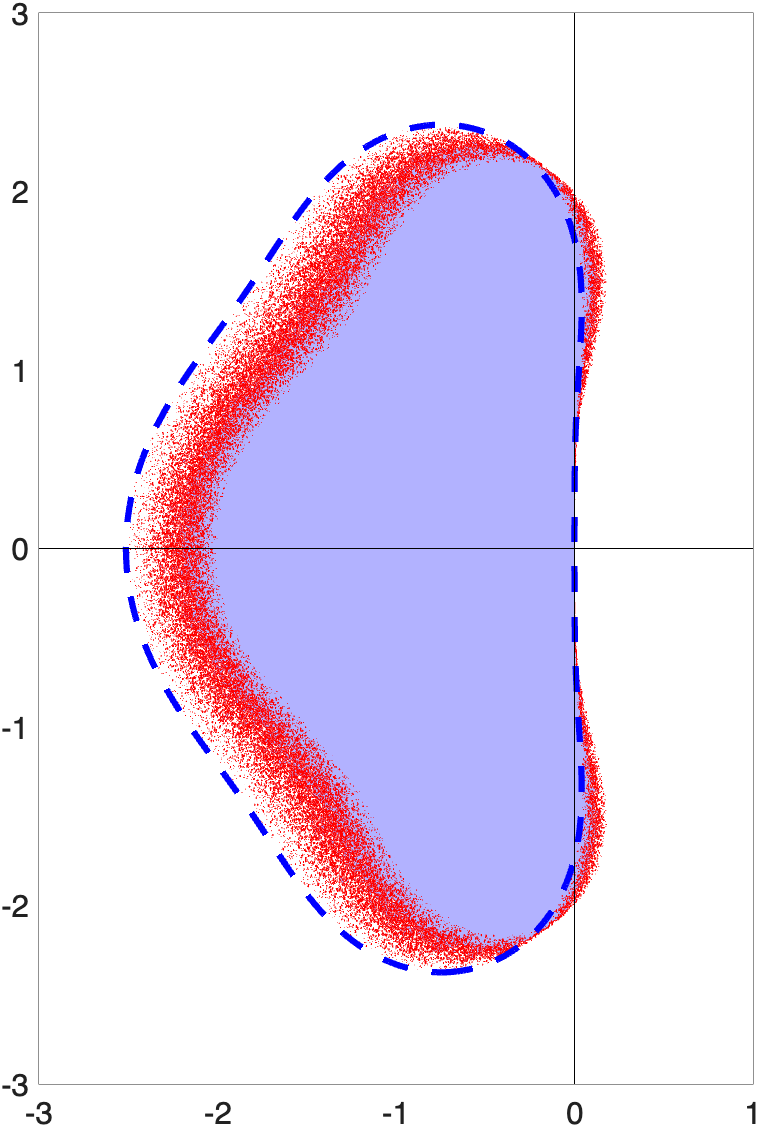}
    \end{subfigure}
\begin{subfigure}[c]{0.27\textwidth}
        \centering
 \includegraphics[width=\textwidth]{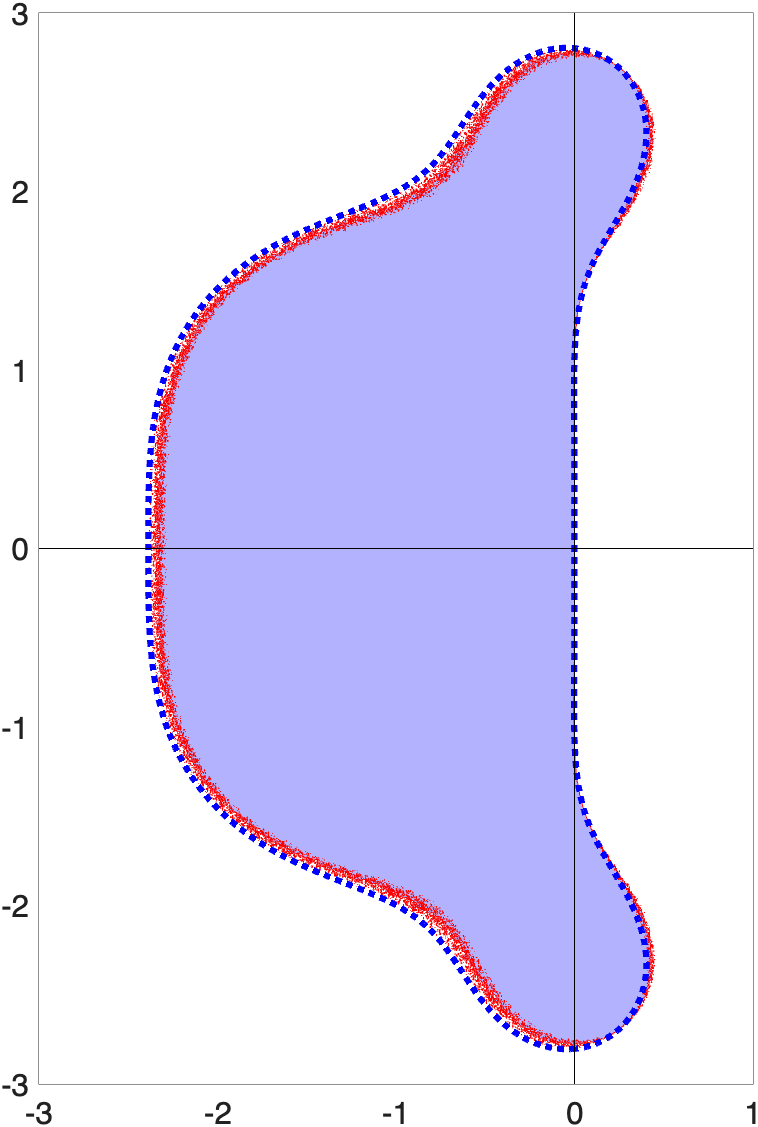}
  \includegraphics[width=\textwidth]{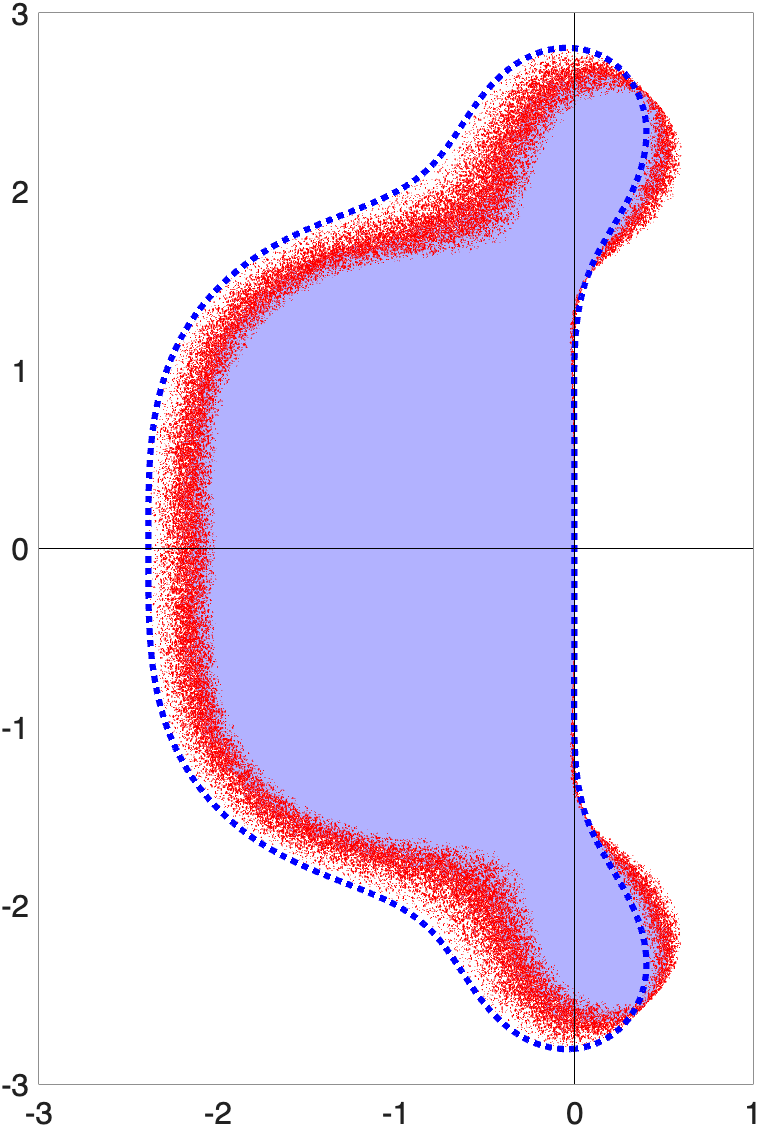}
    \end{subfigure}
    \hfill
     \caption{Perturbed Linear stability regions for the third order MP-TDRK methods. Top: $\ep = 0.1$; 
     Bottom: $\ep = 0.5$. 
     Left: the MP-TDRK2s3p1e method \eqref{MP2s3p1e}.
     Middle: MP-TDRK2s3p2e method \eqref{MP2s3p2e}.
     Right: MP-TDRK3s3p3e method \eqref{MP3s3p3e}.
     }
     \label{fig:PertStab3p}
\end{figure}

Adding the condition \eqref{pert3} we obtain a method with higher 
order $m$ in the perturbation error:

\smallskip

\noindent{\bf TDRK method TDRK3s3p3e:}
This three stage third order method 
\begin{align}  \label{MP3s3p3e}
y^{(1)}    &= u^n + \frac{2}{3} \Delta t F(u^n) 
    + \frac{2}{9} \Delta t^2 \dot{F}_\ep(u^n) \nonumber \\
y^{(2)}    &= u^n + \frac{1}{3} \Delta t F(u^n) 
    + \frac{1}{3} \Delta t F(y^{(1)})  \nonumber \\
u^{n+1} &= u^n + \frac{1}{4} \Delta t F(u^n) 
+ \frac{3}{4} \Delta t F(y^{(2)}) 
\end{align}
will have a final time error of 
\[ E = O(\dt^3) + O(\ep \dt^3).\]

In Figure \ref{fig:stability} (left) we plot the linear stability regions of 
these three third order methods. The stability region of the TDRK2s3p2e method 
\eqref{MP2s3p2e} is smaller than that of the TDRK2s3p1e method 
\eqref{MP2s3p1e}, but it allows for a higher perturbation order $m$.
The TDRK3s3p3e method \eqref{MP3s3p3e} requires an additional stage to obtain $m=3$, 
but this is offset by the larger stability region.
We also want to understand the impact of the 
perturbation error $\ep$ on the linear stability region 
of these methods. To emulate  this, we incorporate an additional
random number of magnitude $(0,\ep)$ into the second derivative
term in the stability polynomial and plot the resulting linear stability regions. These are represented in Figure
\ref{fig:PertStab3p} for $\ep=0.1$ and $\ep = 0.5$, 
where $\ep$  depends upon the precision,  the size
of the problem, and other implementation-dependent factors.
We observe that the larger the perturbation error $\ep$, the smaller 
the linear stability region.

\subsection{Fourth Order Methods} \label{sec:MPTDRK4p}
The  unique explicit two stage fourth order TDRK method  \cite{tsai2010} forms the basis for the first fourth order mixed precision method:\\
\noindent{\bf TDRK method TDRK2s4p1e:}
This two-stage fourth order method 
\begin{align}  \label{MP2s4p1e}
    y^{(1)}    &= u^n + \frac{\Delta t}{2} F(u^n) 
    + \frac{\Delta t^2}{8} \dot{F}_\ep(u^n) \nonumber \\
    u^{n+1} &= u^n + \Delta t F(u^n) + \frac{\Delta t^2}{6}( \dot{F}_\ep(u^n)+
    2\dot{F}_\ep(y^{(1)})),
\end{align}
will have a final time error of 
\[ E = O(\dt^4) + O(\ep \dt).\]

\begin{figure}[tb] 
    \centering
    \begin{subfigure}[c]{0.25\textwidth}
        \centering
\includegraphics[width=\textwidth]{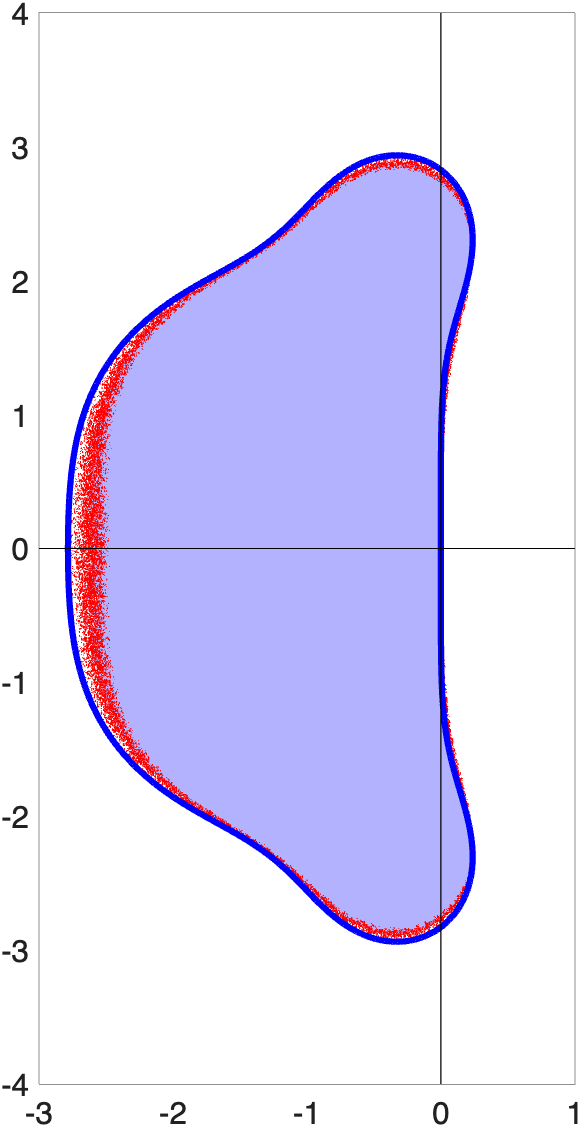}
\includegraphics[width=\textwidth]{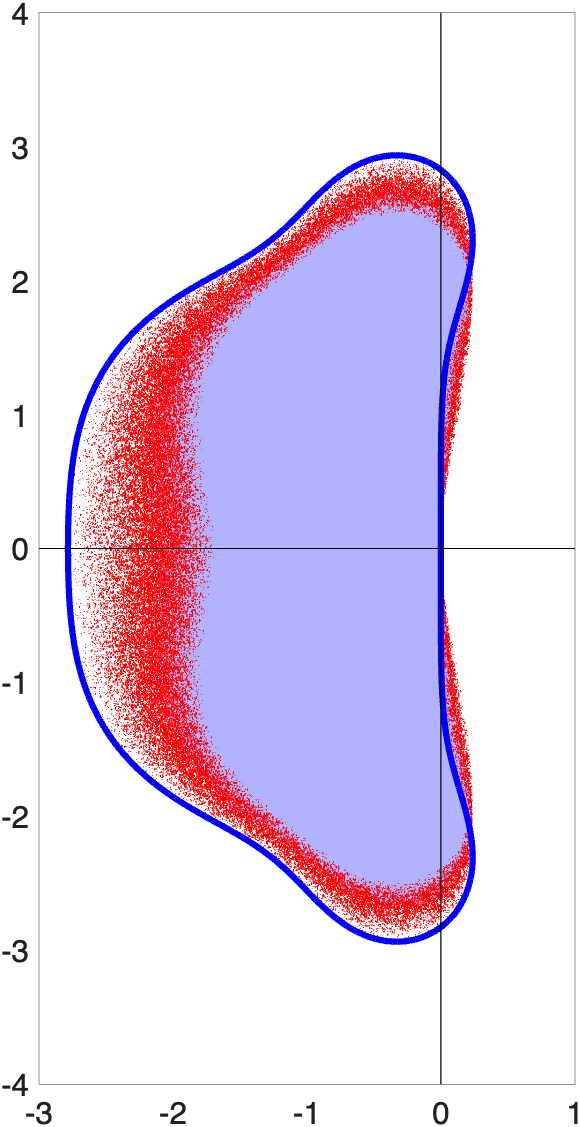}
    \end{subfigure} 
    \hspace{0.25in}
          \begin{subfigure}[c]{0.25\textwidth}
        \centering
 \includegraphics[width=\textwidth]{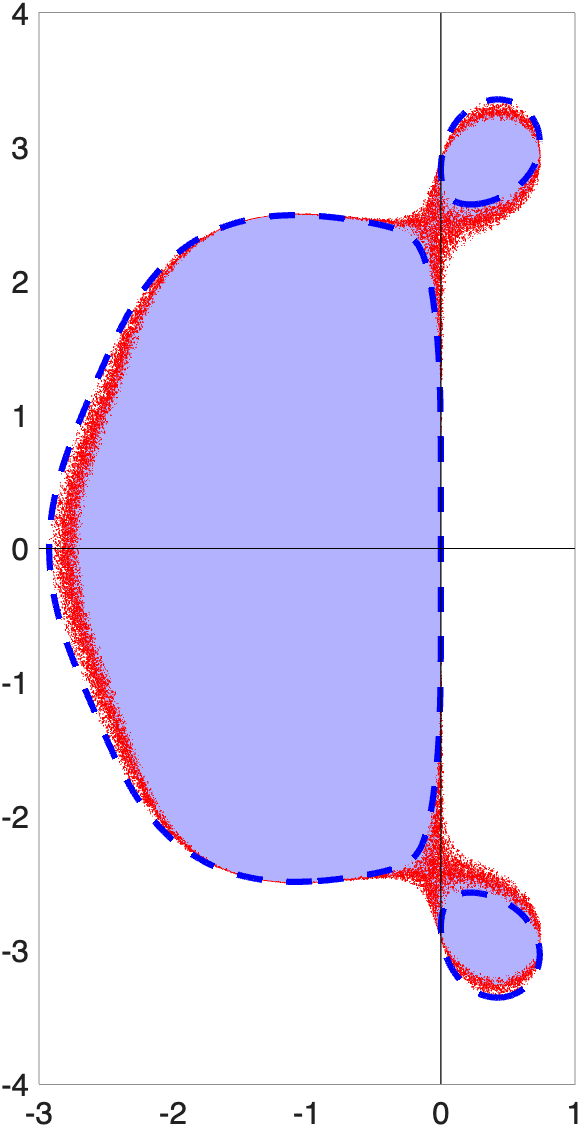}
  \includegraphics[width=\textwidth]{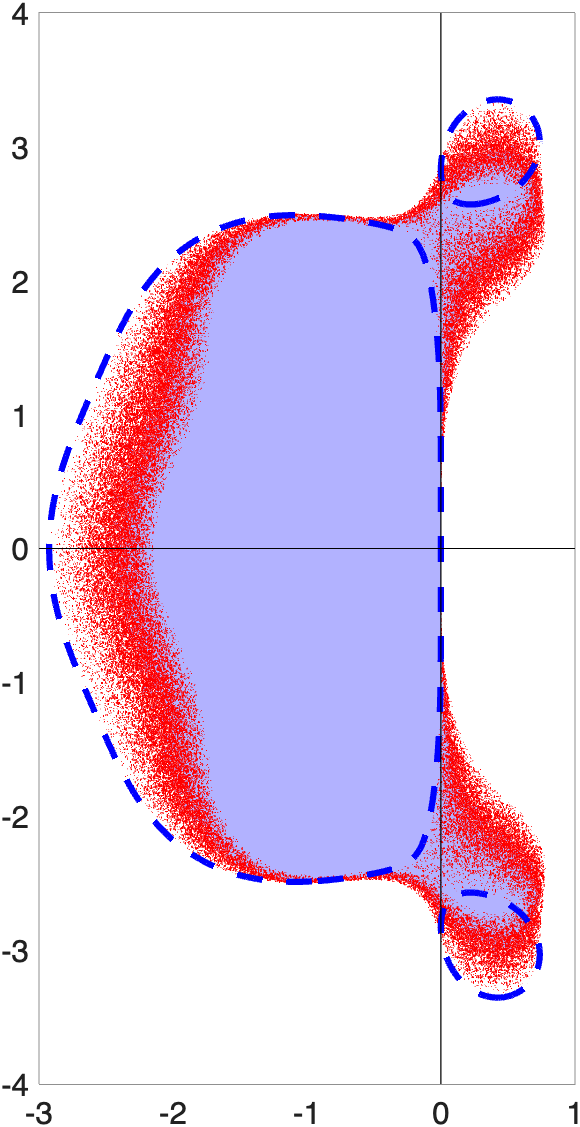}
    \end{subfigure}
     \caption{Perturbed Linear stability regions for the 
     fourth order MP-TDRK methods. Left: two-stage fourth-order
     method \eqref{MP2s4p1e}.
     Right: three-stage fourth-order
     method  \eqref{MP3s4p2e}.
     Top: $\ep = 0.1$; Bottom: $\ep = 0.5$.
     }
     \label{fig:PertStab4p}
\end{figure}

To obtain $m=2$ we require  an additional stage to satisfy 
the additional  perturbation
conditions \eqref{pert2}:\\
\noindent{\bf TDRK method TDRK3s4p2e:}
The  three-stage fourth order method 
\begin{align}  \label{MP3s4p2e}
y^{(1)} &= u^n + \frac{\dt}{2} F(u^n) + \frac{\dt^2}{8} \dot{F}_\ep(u^n)  \nonumber \\
y^{(2)} &= u^n + \dt F(u^n) 
+ \frac{\dt^2}{2} \dot{F}_\ep( y^{(1)} )  \nonumber \\   
u^{n+1} &= u^n + \frac{1}{6} \dt \left( F(u^n) +
    4 F(y^{(1)}) + F(y^{(2)}) \right) ,
\end{align}
will have a final time error of 
\[ E = O(\dt^4) + O(\ep \dt^2).\]
We will need to  keep $ \varepsilon << \dt^2$ to avoid seeing 
an order reduction to second order.  
The price  we pay is additional function evaluations, 
 which contributes to higher computational cost.
Figure \ref{fig:stability} (center) shows that there is no
significant increase in the size of the stability region to offset this additional cost.
Figure \ref{fig:PertStab4p} shows that the $m=1$
TDRK2s4p1e method \eqref{MP2s4p1e} on the left is slightly more impacted by the size of the perturbation than the $m=2$
TDRK3s4p2e method \eqref{MP3s4p2e} on the right.

\subsection{Fifth Order Method} \label{sec:MPTDRK5p}
For the fifth order method we use another of Tsai's methods from \cite{tsai2010}, and implement only the second derivatives in low precision to obtain the mixed precision version. 

\smallskip

\begin{figure}[tb] 
    \centering
\includegraphics[width=0.4\textwidth]{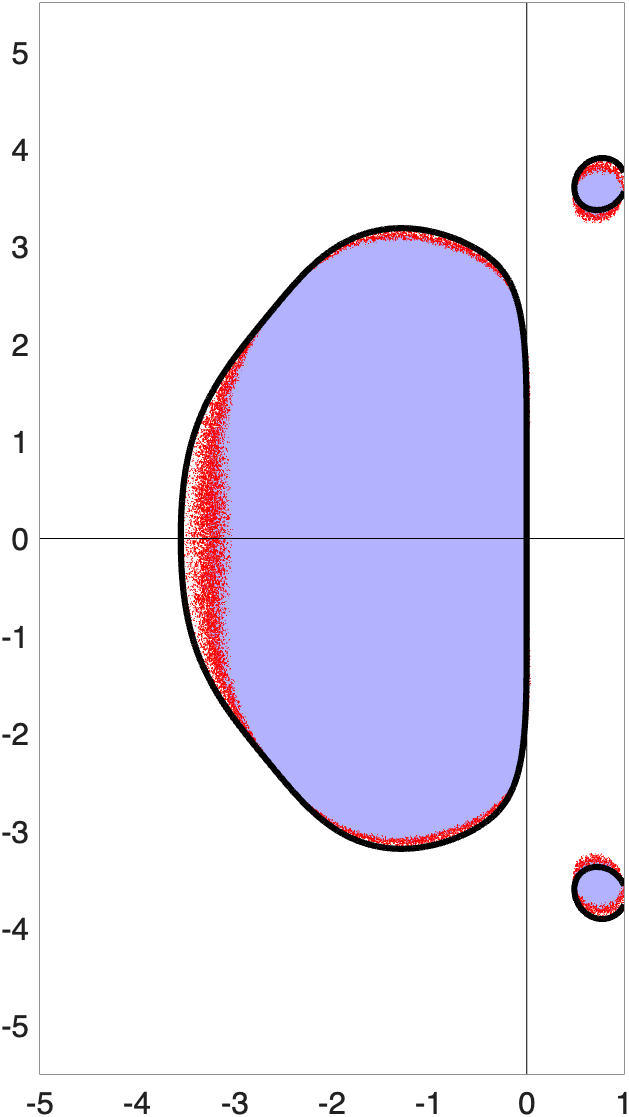}
\hspace{0.25in}
\includegraphics[width=0.4\textwidth]{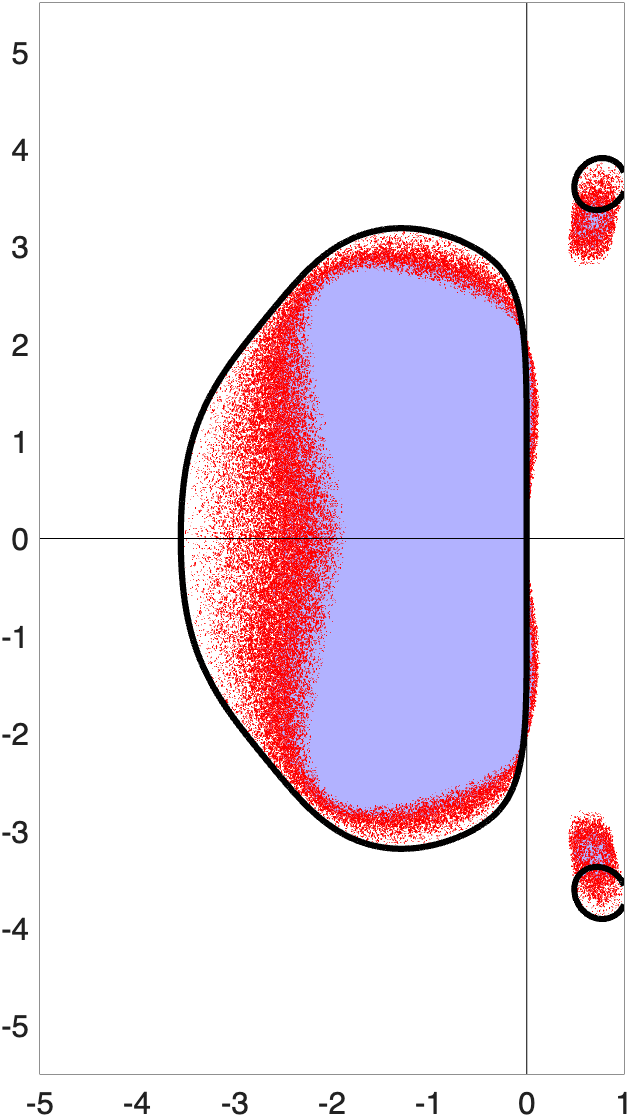}
     \caption{Perturbed Linear stability regions for the 
     fifth order MP-TDRK method  \eqref{Tsai3s5p1e} with $\ep = 0.1$ (left) 
     and $\ep = 0.5$ (right);
     }
     \label{fig:PertStab5p}
\end{figure}

\noindent{\bf TDRK method TDRK3s5p1e:}
The three-stage fifth order method 
\begin{eqnarray} \label{Tsai3s5p1e}
    y^{(1)}     &= & u^n +  \frac{1}{3} \Delta t F(u^n) +  
    \frac{1}{18}\Delta t^2  \dot{F}_\ep(u^n) \nonumber \\
     y^{(2)} &=& u^n +  \frac{4}{5} \Delta t F(u^n)  
    - \frac{2}{125}\Delta t^2  \dot{F}_\ep(u^n) 
    + \frac{42}{125}\Delta t^2  \dot{F}_\ep(y^{(1)}) \nonumber \\
    u^{n+1} &=& u^n +  \Delta t  F(u^n)  
    +  \frac{5}{48} \Delta t^2  \dot{F}_\ep(u^n) +
    \frac{9}{28} \Delta t^2  \dot{F}_\ep(y^{(1)}) \\ 
    &+&
    \frac{25}{336} \Delta t^2 \dot{F}_\ep(y^{(2)}) , \nonumber
\end{eqnarray}
will have a final time error of 
$ E = O(\dt^5) + O(\ep \dt).$
For this fifth order method we require 
$ \ep << \dt^4$ to obtain high order convergence, 
and once $\dt$ is small enough 
we expect to see the first order convergence.

Figure \ref{fig:stability} (right) shows that the linear stabilty region of this method is larger than that 
of the third and fourth order methods. Figure
\ref{fig:PertStab5p} shows that the stability of the  method is very sensitive to the perturbation error $\ep$.

\subsection{Sixth Order Method}  \label{sec:MPTDRK6p}
For the sixth order method we use  Tsai's methods from \cite{tsai2010}. 

\noindent{\bf TDRK method TDRK4s6p1e:}
The four-stage sixth order method 
\begin{eqnarray} \label{MP4s6p1e}
y^{(1)} &=& 
u^n +  \frac{1}{4} \Delta t F(u^n)  + \frac{1}{32} \Delta t^2 \dot{F}_\varepsilon(u^n)  \nonumber \\
y^{(2)} &=& u^n +  \frac{2}{3} \Delta t F(u^n)  
- \frac{2}{81} \dt^2 \dot{F}_\varepsilon(u^n)  + \frac{20}{81} \dt^2 \dot{F}_\varepsilon(y^{(1)})  \nonumber\\ 
y^{(3)} &=& u^n +   \Delta t F(u^n)   +  \frac{5}{4}\dt^2 \dot{F}_\varepsilon(u^n)   \nonumber \\
&& - \frac{6}{5} \dt^2 \dot{F}_\varepsilon(y^{(1)})  
+ \frac{9}{20}\dt^2 \dot{F}_\varepsilon(y^{(2)}) 
\nonumber \\ 
\; \; \; \; \; \; u^{n+1} &=& u^n +  \Delta t  F(u^n)  +  
\frac{3}{40} \Delta t^2 \dot{F}_\varepsilon(u^n)
+  \frac{64}{225} \Delta t^2 \dot{F}_\varepsilon(y^{(1)} )   \\
&& +  \frac{27}{200} \Delta t^2 \dot{F}_\varepsilon(y^{(2)})
+ \frac{1}{180}  \Delta t^2  \dot{F}_\varepsilon(y^{(3)}) . \nonumber
\end{eqnarray}            

This method will have a final time error of 
\[ E = O(\dt^6) + O(\ep \dt),\]
so we expect to see high order when 
$\varepsilon << \dt^5 $, but  when $\dt$ 
is larger  we will see linear convergence. 

\begin{figure}[tb] 
    \centering
\includegraphics[width=0.3\textwidth]{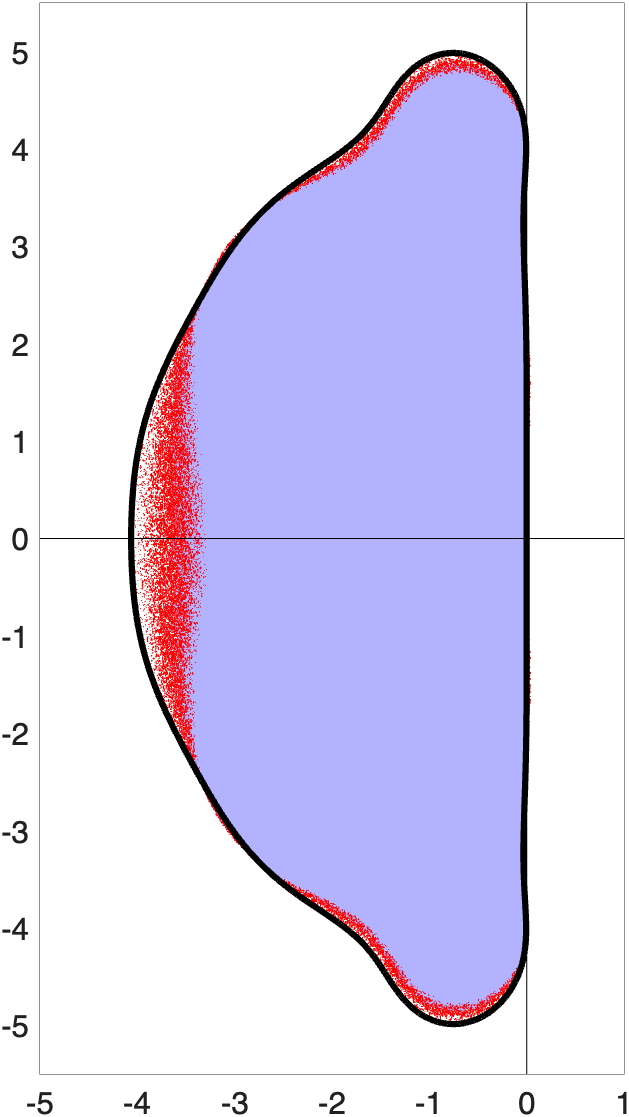}
\hspace{0.25in}
\includegraphics[width=0.3\textwidth]{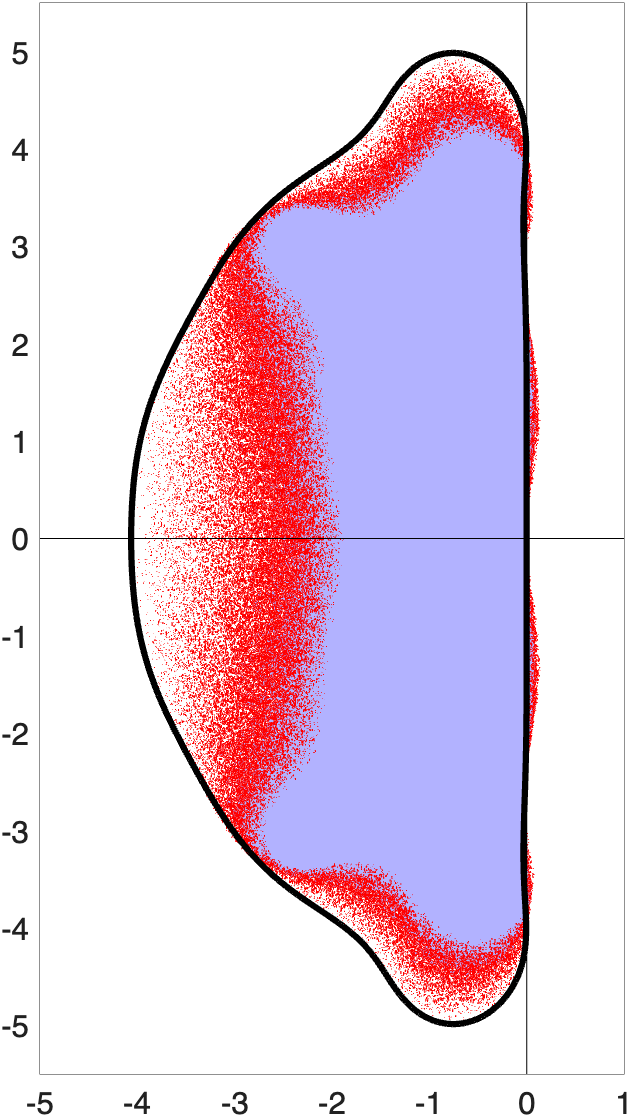}
     \caption{Perturbed Linear stability regions for the 
     sixth order MP-TDRK method  \eqref{MP4s6p1e} with $\ep = 0.1$ (left) 
     and $\ep = 0.5$ (right);
     }
     \label{fig:PertStab6p}
\end{figure}

Figure \ref{fig:stability} (right) shows that the linear stabilty region of this method is larger than that 
all other methods, which may help offset the increased computational cost 
from this four stage method. 
However, Figure \ref{fig:PertStab6p} shows that this method has a 
stability region that is  extremely sensitive to the size of  perturbation error $\ep$.

\section{Numerical Simulations} \label{sec:NumResults}

The focus of this work is the numerical analysis framework described in
Section  \ref{sec:framework}  for developing and analyzing mixed precision TDRK methods.
The  methods studied as well as the new methods designed in Section \ref{sec:MPTDRKmethods} 
are meant to demonstrate that the perturbation order conditions developed 
in Section  \ref{sec:framework} accurately predict the convergence behavior.
These methods are not optimized in any sense,  including efficiency.

The simulations in this section confirm that the impact of 
the low precision perturbations on the accuracy of the solution is as expected
from their design, and that larger perturbations require smaller time-steps 
for stability, as predicted by the perturbed stability regions
in Figures \ref{fig:PertStab3p}, \ref{fig:PertStab4p}, \ref{fig:PertStab5p}, and
\ref{fig:PertStab6p}.
While the use of mixed precision has 
the potential for computational time-savings,
in practice the efficiency of any given methods will depend 
on its design relative to the particular  application and  the relative costs of the operators $F$ and $\dot{F}$.
For these reasons we focus our numerical studies on
the accuracy  and stability of the methods, omitting any 
analysis of their efficiency.  

\subsection{Linear advection with mixed precision}
We consider the linear advection equation
\[ U_t + a U_x = 0 \]
with $a = 1$, and  $x \in [-1, 1]$. We use the  initial condition $U_0(x) = \sin(\pi x)$,
so that the solution is a sine wave at all times, and if we use Fourier spectral methods we have no error from the spatial discretization, so we observe only the time-discretization error. We implement the spatial discretization by defining the Fourier differentiation matrix $\mD_x$ with $N_x$ points in space.

We discretize in space to obtain the system of ODEs
\[ u_t = - a \mD_x u .\] We now use the three different third order  
mixed-precision TDRK methods  \eqref{MP2s3p1e}, \eqref{MP2s3p2e}, and \eqref{MP3s3p3e},
to evolve the solution forward to final time $T_f = 0.5$.

\begin{table}[h!]
\centering
\caption{Errors for mixed precision TDRK2s3p1e}
\resizebox{\textwidth}{!}{
\begin{tabular}{|cc|ccccc|}
\hline
$N_x$ & $\Delta t$ & 64/64 & 64/32 & 32/32 & 64/16 & 16/16 \\
\hline
25 & $10^{-1}$ & $2.04 \times 10^{-3}$ & $1.88 \times 10^{-2}$ & $1.86 \times 10^{-2}$ & $4.83 \times 10^{1}$ & $7.39 \times 10^{1}$ \\
 & $5 \times 10^{-2}$ & $2.54 \times 10^{-4}$ & $2.54 \times 10^{-4}$ & $2.55 \times 10^{-4}$ & $1.11 \times 10^{-2}$ & $3.30 \times 10^{-2}$ \\
 & $2.5 \times 10^{-2}$ & $3.17 \times 10^{-5}$ & $3.18 \times 10^{-5}$ & $3.36 \times 10^{-5}$ & $6.33 \times 10^{-3}$ & $2.91 \times 10^{-2}$ \\
 & $10^{-2}$ & $2.03 \times 10^{-6}$ & $2.15 \times 10^{-6}$ & $4.25 \times 10^{-6}$ & $2.60 \times 10^{-3}$ & $2.91 \times 10^{-2}$ \\
 & $10^{-3}$ & $2.03 \times 10^{-9}$ & $2.37 \times 10^{-8}$ & $3.37 \times 10^{-6}$ & $2.82 \times 10^{-4}$ & $4.33 \times 10^{-2}$ \\
 & $10^{-4}$ & $2.03 \times 10^{-12}$ & $2.56 \times 10^{-9}$ & $1.17 \times 10^{-4}$ & $2.85 \times 10^{-5}$ & $3.86 \times 10^{-1}$ \\
\hline\hline
50 & $10^{-1}$ & $2.24 \times 10^{-3}$ & $2.93 \times 10^{5}$ & $3.18 \times 10^{5}$ & N/A & N/A \\
 & $5 \times 10^{-2}$ & $2.60 \times 10^{-4}$ & $1.08 \times 10^{4}$ & $1.22 \times 10^{4}$ & N/A & N/A \\
 & $2.5 \times 10^{-2}$ & $3.17 \times 10^{-5}$ & $3.25 \times 10^{-5}$ & $4.03 \times 10^{-5}$ & $1.63 \times 10^{-2}$ & $6.69 \times 10^{-2}$ \\
 & $10^{-2}$ & $2.03 \times 10^{-6}$ & $2.98 \times 10^{-6}$ & $1.24 \times 10^{-5}$ & $1.05 \times 10^{-2}$ & $8.25 \times 10^{-2}$ \\
 & $10^{-3}$ & $2.03 \times 10^{-9}$ & $1.14 \times 10^{-7}$ & $9.86 \times 10^{-6}$ & $8.45 \times 10^{-4}$ & $9.30 \times 10^{-2}$ \\
 & $10^{-4}$ & $2.03 \times 10^{-12}$ & $1.18 \times 10^{-8}$ & $1.48 \times 10^{-4}$ & $8.56 \times 10^{-5}$ & $3.87 \times 10^{-1}$ \\
\hline\hline
100 & $10^{-1}$ & $3.13 \times 10^{3}$ & $1.25 \times 10^{12}$ & $1.25 \times 10^{12}$ & N/A & N/A \\
 & $5 \times 10^{-2}$ & $2.63 \times 10^{9}$ & $9.83 \times 10^{17}$ & $9.90 \times 10^{17}$ & N/A & N/A \\
 & $2.5 \times 10^{-2}$ & $1.99 \times 10^{7}$ & $6.24 \times 10^{15}$ & $6.74 \times 10^{15}$ & N/A & N/A \\
 & $10^{-2}$ & $2.03 \times 10^{-6}$ & $8.68 \times 10^{-6}$ & $1.81 \times 10^{-5}$ & $8.05 \times 10^{-2}$ & $1.53 \times 10^{-1}$ \\
 & $10^{-3}$ & $2.03 \times 10^{-9}$ & $6.01 \times 10^{-7}$ & $1.31 \times 10^{-5}$ & $7.45 \times 10^{-3}$ & $1.57 \times 10^{-1}$ \\
 & $10^{-4}$ & $2.06 \times 10^{-12}$ & $6.15 \times 10^{-8}$ & $1.53 \times 10^{-4}$ & $7.37 \times 10^{-4}$ & $2.02 \times 10^{-1}$ \\
\hline
\end{tabular}
}
\label{tab:linadv1e}
\end{table}

\begin{table}[h!]
\centering
\caption{Errors for mixed precision TDRK2s3p2e}
\resizebox{\textwidth}{!}{
\begin{tabular}{|cc|ccccc|}
\hline
$N_x$ & $\Delta t$ & 64/64 & 64/32 & 32/32 & 64/16 & 16/16 \\
\hline
25 & $10^{-1}$ & $2.03 \times 10^{-3}$ & $5.42 \times 10^{-3}$ & $4.74 \times 10^{-3}$ & $1.07 \times 10^{1}$ & $1.18 \times 10^{1}$ \\
 & $5 \times 10^{-2}$ & $2.54 \times 10^{-4}$ & $2.54 \times 10^{-4}$ & $2.52 \times 10^{-4}$ & $4.63 \times 10^{-3}$ & $3.35 \times 10^{-2}$ \\
 & $2.5 \times 10^{-2}$ & $3.16 \times 10^{-5}$ & $3.17 \times 10^{-5}$ & $3.15 \times 10^{-5}$ & $9.44 \times 10^{-4}$ & $2.84 \times 10^{-2}$ \\
 & $10^{-2}$ & $2.03 \times 10^{-6}$ & $2.04 \times 10^{-6}$ & $3.55 \times 10^{-6}$ & $1.54 \times 10^{-4}$ & $3.13 \times 10^{-2}$ \\
 & $10^{-3}$ & $2.03 \times 10^{-9}$ & $2.15 \times 10^{-9}$ & $3.21 \times 10^{-6}$ & $1.58 \times 10^{-6}$ & $8.20 \times 10^{-2}$ \\
 & $10^{-4}$ & $2.03 \times 10^{-12}$ & $3.40 \times 10^{-12}$ & $8.77 \times 10^{-6}$ & $1.58 \times 10^{-8}$ & $6.26 \times 10^{-1}$ \\
\hline\hline
50 & $10^{-1}$ & $2.03 \times 10^{-3}$ & $7.47 \times 10^{2}$ & $7.04 \times 10^{2}$ & N/A & N/A \\
 & $5 \times 10^{-2}$ & $2.54 \times 10^{-4}$ & $2.73 \times 10^{2}$ & $2.74 \times 10^{2}$ & N/A & N/A \\
 & $2.5 \times 10^{-2}$ & $3.16 \times 10^{-5}$ & $3.28 \times 10^{-5}$ & $3.51 \times 10^{-5}$ & $7.54 \times 10^{-3}$ & $7.13 \times 10^{-2}$ \\
 & $10^{-2}$ & $2.03 \times 10^{-6}$ & $2.17 \times 10^{-6}$ & $9.83 \times 10^{-6}$ & $1.56 \times 10^{-3}$ & $7.71 \times 10^{-2}$ \\
 & $10^{-3}$ & $2.03 \times 10^{-9}$ & $3.85 \times 10^{-9}$ & $9.56 \times 10^{-6}$ & $1.68 \times 10^{-5}$ & $9.46 \times 10^{-2}$ \\
 & $10^{-4}$ & $2.03 \times 10^{-12}$ & $2.00 \times 10^{-11}$ & $1.23 \times 10^{-5}$ & $1.72 \times 10^{-7}$ & $5.66 \times 10^{-1}$ \\
\hline\hline
100 & $10^{-1}$ & $1.36 \times 10^{-1}$ & $5.95 \times 10^{7}$ & $6.01 \times 10^{7}$ & N/A & N/A \\
 & $5 \times 10^{-2}$ & $5.36 \times 10^{3}$ & $2.58 \times 10^{12}$ & $2.55 \times 10^{12}$ & N/A & N/A \\
 & $2.5 \times 10^{-2}$ & $3.27 \times 10^{3}$ & $1.40 \times 10^{12}$ & $1.55 \times 10^{12}$ & N/A & N/A \\
 & $10^{-2}$ & $2.03 \times 10^{-6}$ & $2.92 \times 10^{-6}$ & $1.34 \times 10^{-5}$ & $9.32 \times 10^{-3}$ & $9.76 \times 10^{-2}$ \\
 & $10^{-3}$ & $2.03 \times 10^{-9}$ & $1.74 \times 10^{-8}$ & $1.35 \times 10^{-5}$ & $1.47 \times 10^{-4}$ & $1.06 \times 10^{-1}$ \\
 & $10^{-4}$ & $2.05 \times 10^{-12}$ & $1.62 \times 10^{-10}$ & $1.35 \times 10^{-5}$ & $1.48 \times 10^{-6}$ & $4.19 \times 10^{-1}$ \\
\hline
\end{tabular}
}
\label{tab:linadv2e}
\end{table}

\begin{table}[h!]
\centering
\caption{Errors for mixed precision TDRK3s3p}
\resizebox{\textwidth}{!}{
\begin{tabular}{|cc|ccccc|}
\hline
$N_x$ & $\Delta t$ & 64/64 & 64/32 & 32/32 & 64/16 & 16/16 \\
\hline
25 & $10^{-1}$ & $6.95 \times 10^{-4}$ & $2.69 \times 10^{-3}$ & $2.13 \times 10^{-3}$ & $1.10 \times 10^{1}$ & $3.37 \times 10^{0}$ \\
 & $5 \times 10^{-2}$ & $8.51 \times 10^{-5}$ & $8.53 \times 10^{-5}$ & $8.61 \times 10^{-5}$ & $1.68 \times 10^{-3}$ & $2.82 \times 10^{-2}$ \\
 & $2.5 \times 10^{-2}$ & $1.06 \times 10^{-5}$ & $1.06 \times 10^{-5}$ & $1.24 \times 10^{-5}$ & $2.69 \times 10^{-4}$ & $2.91 \times 10^{-2}$ \\
 & $10^{-2}$ & $6.76 \times 10^{-7}$ & $6.76 \times 10^{-7}$ & $3.25 \times 10^{-6}$ & $1.81 \times 10^{-5}$ & $3.08 \times 10^{-2}$ \\
 & $10^{-3}$ & $6.76 \times 10^{-10}$ & $6.76 \times 10^{-10}$ & $3.40 \times 10^{-6}$ & $1.81 \times 10^{-8}$ & $3.74 \times 10^{-2}$ \\
 & $10^{-4}$ & $6.76 \times 10^{-13}$ & $6.77 \times 10^{-13}$ & $1.01 \times 10^{-5}$ & $1.80 \times 10^{-11}$ & $3.98 \times 10^{-1}$ \\
\hline\hline
50 & $10^{-1}$ & $7.28 \times 10^{-4}$ & $4.52 \times 10^{4}$ & $4.28 \times 10^{4}$ & N/A & N/A \\
 & $5 \times 10^{-2}$ & $8.52 \times 10^{-5}$ & $1.47 \times 10^{2}$ & $8.42 \times 10^{1}$ & N/A & N/A \\
 & $2.5 \times 10^{-2}$ & $1.06 \times 10^{-5}$ & $1.10 \times 10^{-5}$ & $1.84 \times 10^{-5}$ & $2.80 \times 10^{-3}$ & $6.74 \times 10^{-2}$ \\
 & $10^{-2}$ & $6.76 \times 10^{-7}$ & $7.15 \times 10^{-7}$ & $1.05 \times 10^{-5}$ & $2.83 \times 10^{-4}$ & $8.01 \times 10^{-2}$ \\
 & $10^{-3}$ & $6.76 \times 10^{-10}$ & $7.17 \times 10^{-10}$ & $1.03 \times 10^{-5}$ & $2.72 \times 10^{-7}$ & $9.42 \times 10^{-2}$ \\
 & $10^{-4}$ & $6.72 \times 10^{-13}$ & $7.26 \times 10^{-13}$ & $2.07 \times 10^{-5}$ & $2.82 \times 10^{-10}$ & $4.92 \times 10^{-1}$ \\
\hline\hline
100 & $10^{-1}$ & $3.81 \times 10^{2}$ & $1.63 \times 10^{11}$ & $1.54 \times 10^{11}$ & N/A & N/A \\
 & $5 \times 10^{-2}$ & $3.18 \times 10^{7}$ & $1.69 \times 10^{16}$ & $1.59 \times 10^{16}$ & N/A & N/A \\
 & $2.5 \times 10^{-2}$ & $1.15 \times 10^{3}$ & $6.86 \times 10^{11}$ & $3.87 \times 10^{11}$ & N/A & N/A \\
 & $10^{-2}$ & $6.77 \times 10^{-7}$ & $1.05 \times 10^{-6}$ & $1.39 \times 10^{-5}$ & $4.29 \times 10^{-3}$ & $1.02 \times 10^{-1}$ \\
 & $10^{-3}$ & $6.76 \times 10^{-10}$ & $1.41 \times 10^{-9}$ & $1.34 \times 10^{-5}$ & $5.70 \times 10^{-6}$ & $1.10 \times 10^{-1}$ \\
 & $10^{-4}$ & $7.09 \times 10^{-13}$ & $1.43 \times 10^{-12}$ & $2.04 \times 10^{-5}$ & $5.52 \times 10^{-9}$ & $1.92 \times 10^{-1}$ \\
\hline
\end{tabular}
}
\label{tab:linadv3e}
\end{table}

Tables \ref{tab:linadv1e}, \ref{tab:linadv2e}, and \ref{tab:linadv3e}
show the final time maximum norm errors 
compared to the exact  solution  $U(x,t) = \sin(\pi (x-a t)) $
for 
$N_x = 25,50,100$, for each of the methods
\eqref{MP2s3p1e}, \eqref{MP2s3p2e}, and \eqref{MP3s3p3e}, respectively.
These results confirm that the methods display 
the perturbation order as predicted by the theory.
This results in smaller errors
from the mixed precision implementation of the methods
with perturbation errors $\ep \dt^m$ as $m$ increases.
 It is notable that while 
for  $N_x=25$ the mixed precision 
method converges for all values of $\dt$, 
as $N_x$ gets larger we require smaller values of $\dt$ for convergence. This effect is more pronounced as we move to lower precision. The mixed method inherits the stability properties of the 
low-precision method, so that while mixed precision improves the magnitude of the errors, it does not improve the value of $\dt$ needed for stability of the method. However, the mixed precision method ensures that the method will converge as $\dt$ 
continues to gets much smaller, which is not the case for the lower precision method.

This behavior explained by the stability regions plotted in Figure \ref{fig:PertStab3p}. Two factors are at play here.
First, as the perturbation grows larger
(i.e. the precision is lower), the linear stability region shrinks.
In addition, the value of $\ep$ is not only related to the 
precision, but also to the size of the problem, as we explained in
Remark \ref{rmk2}. For this reason, even for the same precision, $\ep$ is larger as $N_x$ is larger, resulting in smaller stability regions.

\begin{figure}[t] 
    \centering
    \begin{subfigure}[b]{0.46\textwidth}
        \centering
        \includegraphics[width=\textwidth]{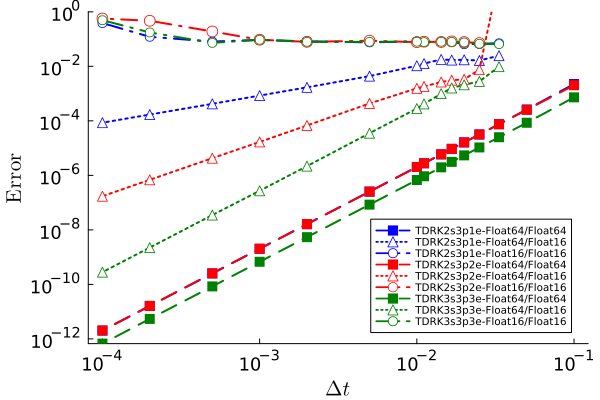}
\includegraphics[width=\textwidth]{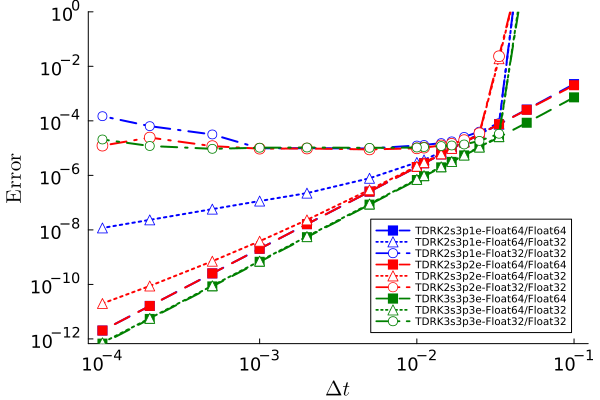}
    \end{subfigure}
    \hfill
           \begin{subfigure}[b]{0.46\textwidth}
        \centering
        \includegraphics[width=\textwidth]{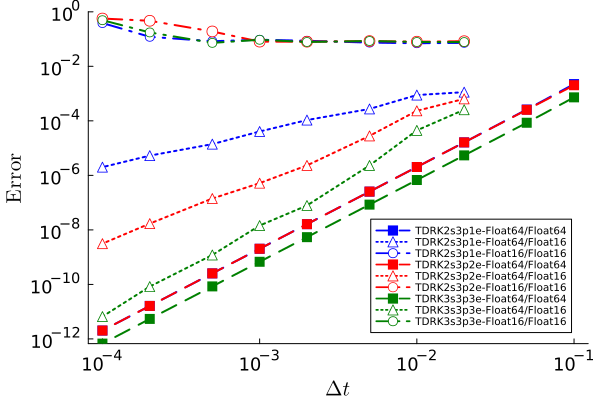}
   \includegraphics[width=\textwidth]{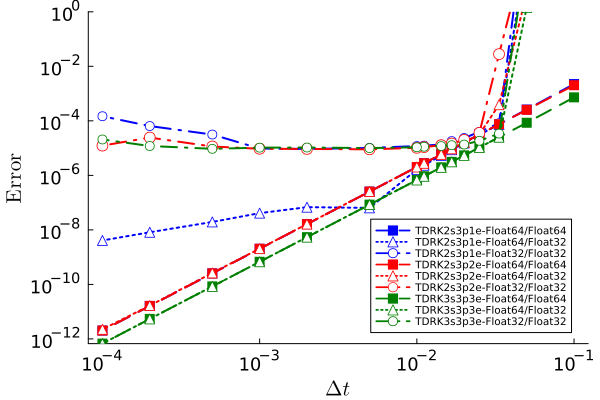}
    \end{subfigure}
    \hfill
     \caption{Mixed double/half precision third order methods applied to the linear advection equation
     for $N_x = 50$. On the left is \textbf{Implementation 1} and on the right is the more accurate but also more costly \textbf{Implementation 2}. Top: double/half;
     Bottom: double/single. Note that NaNs and INFs do not plot, so that the lines end at smaller $\dt$.}
     \label{fig:linadv3p}
\end{figure}

Figure \ref{fig:linadv3p} 
shows the final time maximum-norm errors compared to the exact  solution for $U(x,t) = \sin(\pi (x-a t)) $
for $N_x = 50$ with double/half mixed precision  (left) and 
with double/single mixed precision (right). 
Note that as $\dt$ gets very small,
the low precision methods do not converge,
and indeed the error grows. The mixed precision methods converge
as $\ep \dt^m$ for $m=1,2,3$ for \eqref{MP2s3p1e}, \eqref{MP2s3p2e}, and \eqref{MP3s3p3e}, respectively.
The top figures show double precision mixed with half 
precision, and the bottom figures show double precision 
mixed with single  precision. We observe that for the 
mixed double/single method the errors from the $m=2$ method 
\eqref{MP2s3p2e} look very close to those from
the high precision method,  and the error from the $m=3$ method \eqref{MP3s3p3e} are indistinguishable  from those from the high precision method.

The two sides  of Figure \ref{fig:linadv3p} show the 
impact of different mixed precision implementation.
On the left we have {\bf Implementation 1}, in which the 
compute the second derivative matrix in low precision by 
using low precision FFTs. This is computationally faster
but results in a buildup of
the low precision rounding errors and a larger $\ep$.
This is the primary implementation we use in this work.
On the right we show the result of a more accurate 
{\bf Implementation 2}, in which we
compute the second derivative matrix in high precision initially,
and then downcast it to lower precision and use it in the 
computation of $\dot{F}$. We observe that the mixed precision error lines in {\bf Implementation 2} are shifted downwards by one to two
orders of magnitude compared to {\bf Implementation 1}. This 
shows that the smaller $\ep$ results in significantly smaller 
errors.

\begin{figure}[b]
    \centering
    \begin{subfigure}[b]{0.45\textwidth}
        \centering
        \includegraphics[width=\textwidth]{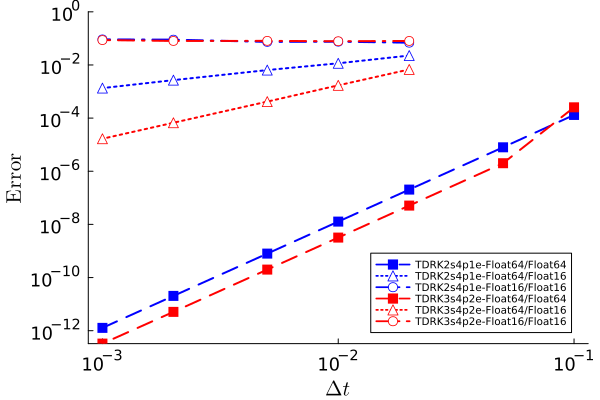}
    \end{subfigure}
    \hfill
    \begin{subfigure}[b]{0.45\textwidth}
        \centering
        \includegraphics[width=\textwidth]{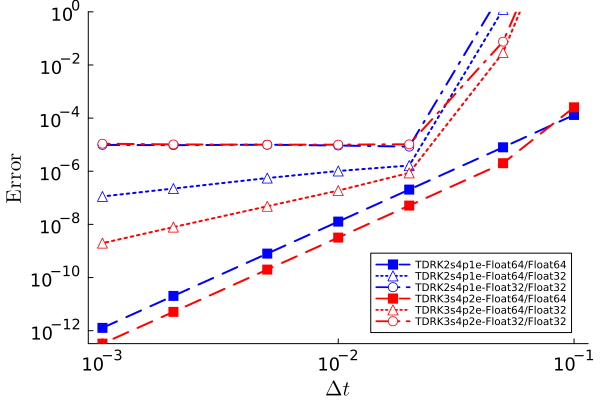}
    \end{subfigure}
    \hfill
    \caption{Mixed precision fourth order methods applied to the linear advection equation
     for $N_x = 50$. On the left is double/half and in the right is double/single.}\label{fig:linadv4p}
\end{figure}

Figure \ref{fig:linadv4p} shows the error plots for 
$N_x =50$ for the 
the two fourth order methods \eqref{MP2s4p1e} and \eqref{MP3s4p2e}
that have perturbation orders $O(\ep \dt)$ and $O(\ep \dt^2)$, respectively.
On the left we see the behavior with 
double precision  mixed with half precision, and on the right the 
double precision  mixed with single precision. 
The errors are significantly lower in the double/single implementation, and we see a larger $\dt$ is allowed for stability.
We observe that, as predicted by the theory,  the TDRK2s4p1e method \eqref{MP2s4p1e}
has first order convergence once
$\dt$ is small enough compared to $\ep$, while 
the TDRK3s4p2e method
\eqref{MP3s4p2e}  has second order convergence.

\begin{figure}[t]
    \centering
    \begin{subfigure}[b]{0.4\textwidth}
        \centering
        \includegraphics[width=\textwidth]{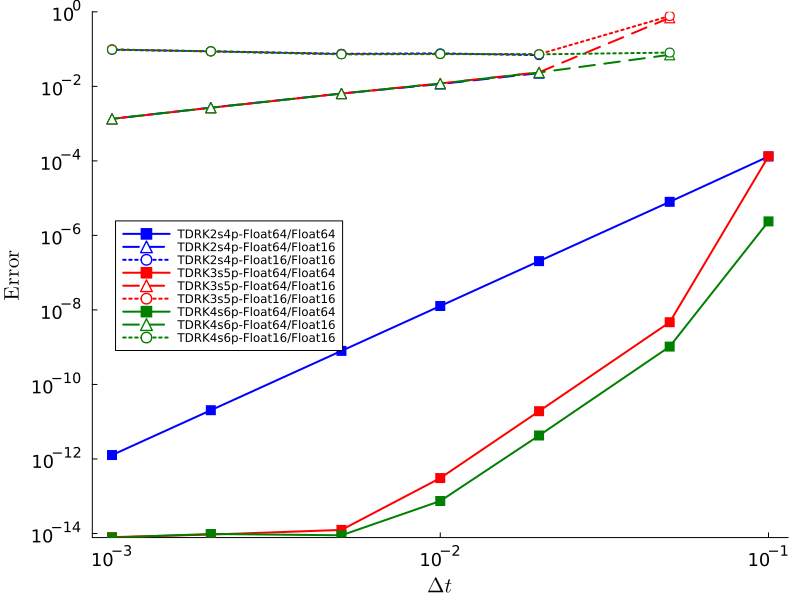}
    \end{subfigure}
    \hspace{10pt}
    \begin{subfigure}[b]{0.4\textwidth}
        \centering
        \includegraphics[width=\textwidth]{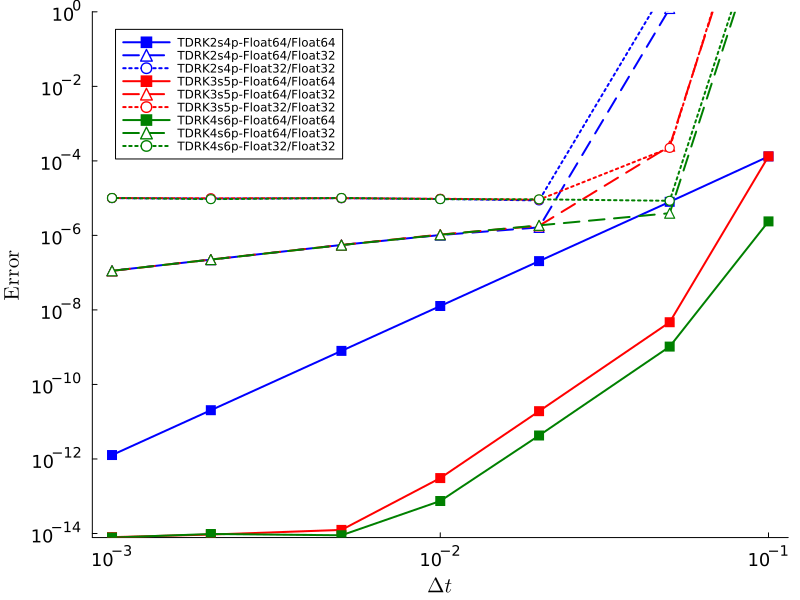}
    \end{subfigure}
    \hfill
    \caption{Mixed precision  methods of order $p=4,5,6$ and perturbation order $m=1$
    applied to the linear advection equation for  $N_x = 50$. Left has mixed double/half precision and right has mixed double/single precision.}\label{fig:TDRK4p5p6p1e}
\end{figure}

Figure \ref{fig:TDRK4p5p6p1e} shows the error plots for $N_x =50$ for the 
the fourth, fifth, and sixth order mixed precision TDRK 
methods \eqref{MP2s4p1e}, \eqref{Tsai3s5p1e}, \eqref{MP4s6p1e}, 
that have perturbation errors $O(\ep \dt)$. On the left we see the behavior with 
double precision  mixed with half precision, and on the right the 
double precision  mixed with single precision. 
These methods all converge linearly in the mixed precision 
setting.
The figure confirms the design order accuracy of $p=4,6$.
For $p=5$ we have higher order convergence than we expect.
Closer analysis of the coefficients in \eqref{Tsai3s5p1e}
show that the order condition corresponding to the linear
tree (the tall tree) is satisfied. This is due to the simplified
structure inherent in Tsai's methods, which causes several of the sixth order trees to be satisfied, as well, due to the high stage order that is enforced in Tsai's methods. 
Thus, the $p=5$ method
is in fact sixth order when applied to linear problems.
We also observe first order convergence in the mixed precision methods, as expected, with the double/single errors four orders
of magnitude smaller than the double/half errors.

\subsection{Inviscid Burgers' equation with Fourier 
spectral methods}

In this subsection we solve the inviscid 
Burgers' equation 
\begin{eqnarray}
    U_t  + \left(\frac{1}{2} U^2 \right)_x = 0,
\end{eqnarray}
for $x \in [-1, 1]$, with the initial condition 
$U_0(x) = \frac{1}{2} + \frac{1}{4}\sin(\pi x)$.

We discretize in space using the spectral differentiation 
matrix described above, to obtain
\begin{eqnarray} \label{BurgersODE}
    u_t  &=& F(u)  = - \mD_x f(u) 
\end{eqnarray}
and compute the second derivative
\[  \dot{F}(u) = \mD_x \left(u F(u)\right).\]

We evolve the solution forward in time to $T_f =0.5$ 
using the different mixed-precision TDRK methods. 
For each $N_x$, the reference solution is computed in quad precision using the Shu-Osher SSP33 method \cite{Shu1988a}
with the same $N_x$ and a refined 
time step of $\Delta t = 10^{-7}$.

\begin{figure}[t]
    \centering
    \begin{subfigure}[b]{0.45\textwidth}
        \centering
        \includegraphics[width=\textwidth]{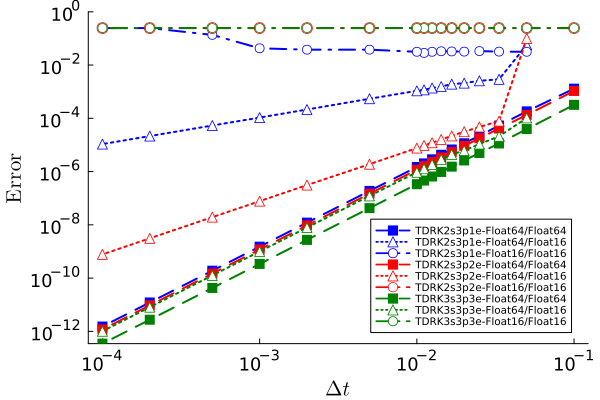}
    \end{subfigure}
    \hfill
    \begin{subfigure}[b]{0.45\textwidth}
        \centering
        \includegraphics[width=\textwidth]{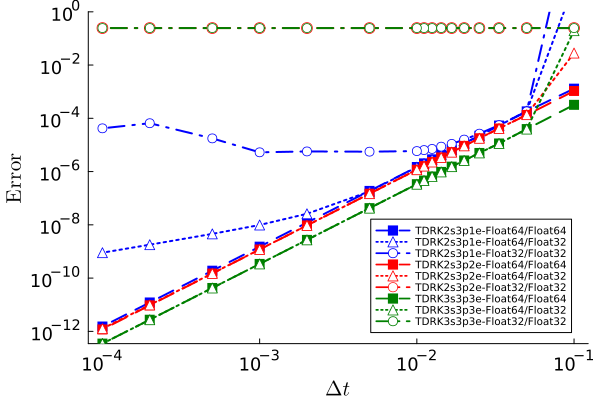}
    \end{subfigure}
    \hfill
    \caption{Mixed precision third order methods 
    applied to Burger's equation  for $N_x = 50$. 
         On the left is double/half and in the right is double/single.
   Method \eqref{MP2s3p1e} in blue,
   method \eqref{MP2s3p2e} in red,
   and method \eqref{MP3s3p3e} in green.
}\label{fig:burgers3p}
\end{figure}

Figure \ref{fig:burgers3p}  shows the 
errors from Burgers' equation evolved with the 
third order  methods  \eqref{MP2s3p1e}, \eqref{MP2s3p2e},
 and \eqref{MP3s3p3e} for double/half (left) 
 double/single (right). This figure shows the impact of the
perturbation order $\ep \dt^m$
with increasing $m$.  Notably, for both 
precisions, the errors from the $O(\ep \dt^3) $
mixed precision method \eqref{MP3s3p3e} 
are indistinguishable from the full precision implementation.
For double/single, the errors from the 
$O(\ep \dt^2) $ mixed precision method \eqref{MP2s3p2e} 
are also indistinguishable from the full precision implementation,
and the errors from the $O(\ep \dt) $ mixed precision method \eqref{MP2s3p1e}
converge at third order until below $\dt < 10^{-2}$, and then they switch to 
first order, as the perturbation error dominates.

\begin{figure}[b]
    \centering
    \begin{subfigure}[b]{0.45\textwidth}
        \centering
    \includegraphics[width=\textwidth]{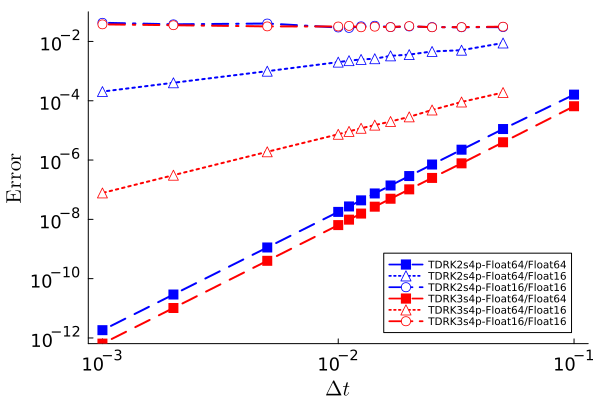}
    \end{subfigure}
    \hfill
    \begin{subfigure}[b]{0.45\textwidth}
        \centering
    \includegraphics[width=\textwidth]{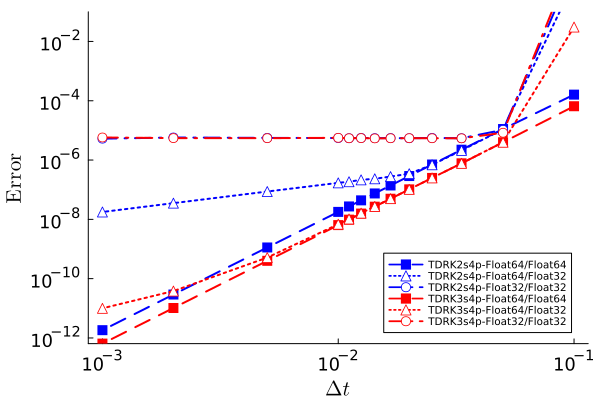}
    \end{subfigure}
    \hfill
    \caption{Errors from Burger's equation for $N_x = 50$ evolved with the 
     mixed precision fourth order methods \eqref{MP2s4p1e} and \eqref{MP3s4p2e}. 
    On the left is double/half and in the right is double/single.} \label{fig:burgers4p}
\end{figure}

\begin{figure}[t]
    \centering
    \begin{subfigure}[b]{0.48\textwidth}
        \centering        \includegraphics[width=\textwidth]{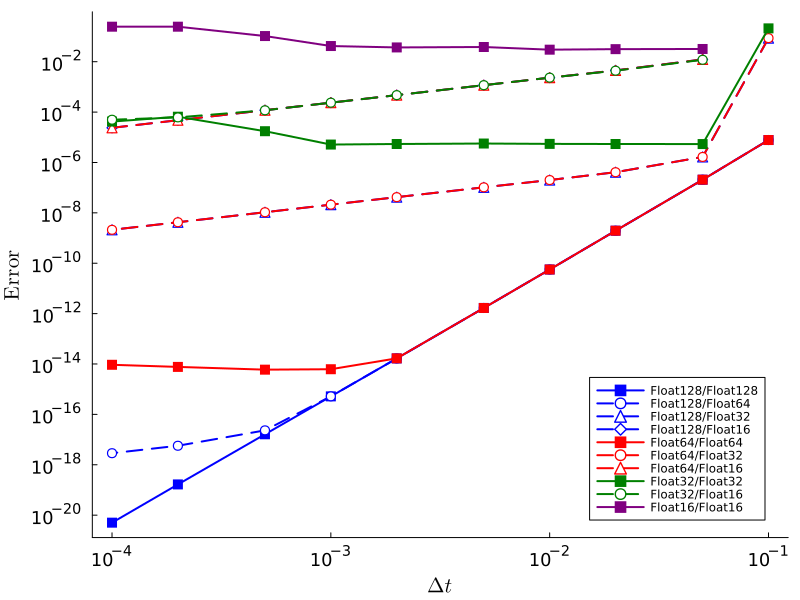}
    \end{subfigure}
    \hfill
    \begin{subfigure}[b]{0.48\textwidth}
        \centering   \includegraphics[width=\textwidth]{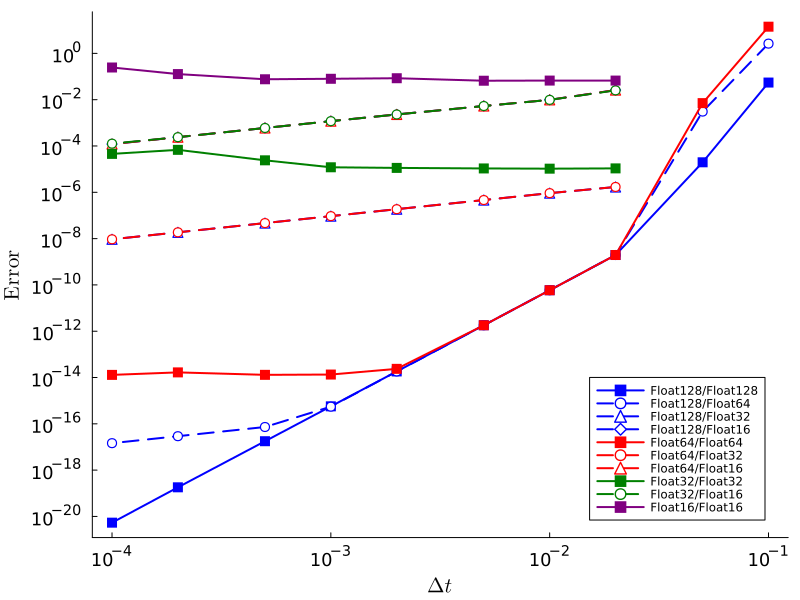}
    \end{subfigure}
    \caption{Mixed precision  fifth order method TDRK3s5p1e \eqref{MP3s4p2e} applied to Burgers' equation for (left) $N_x = 50$, and (right) $N_x = 100$.}\label{fig:burgers5p}
\end{figure}

Figure \ref{fig:burgers4p}  shows the 
errors from Burgers' equation evolved with the 
fourth order methods \eqref{MP2s4p1e} in blue and 
\eqref{MP3s4p2e} in red. 
As expected, the methods are fourth order in full precision, first order in mixed precision, and not convergent in low precision.
While in  mixed
double/half precision (left)
there is a a difference of over four orders of magnitude between the mixed precision and full precision implementation of method 
\eqref{MP3s4p2e} (in red), in 
mixed double/single precision (right) the errors from  
method \eqref{MP3s4p2e} are fifth order until well below $10^{-2}$, at which point they switch over to second order. This is a clear illustration of the errors $E = O(\dt^4) + O(\ep \dt^2) $, where the fourth order error dominates until $\dt$ becomes small enough compared to $\ep$ that the second order error becomes dominant. 

\begin{figure}[b]
    \centering
    \begin{subfigure}[b]{0.48\textwidth}
        \centering        \includegraphics[width=\textwidth]{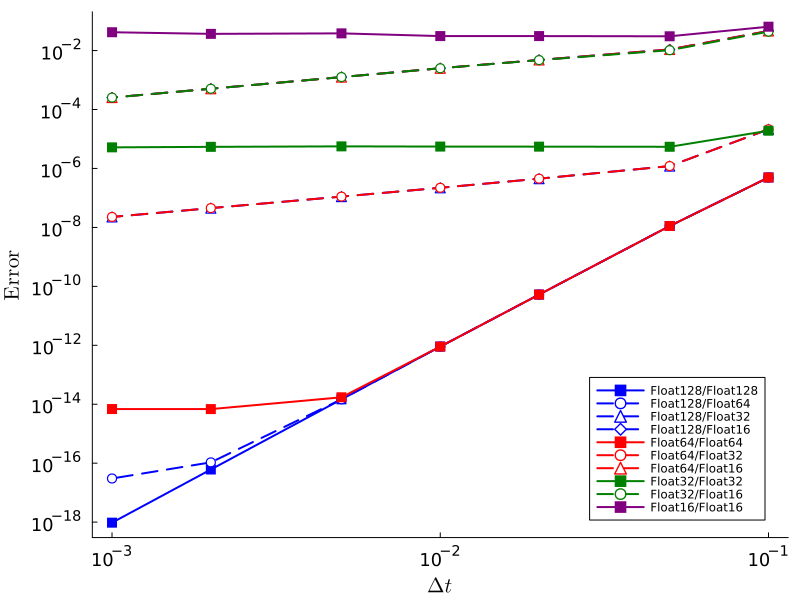}
    \end{subfigure}
    \hfill
    \begin{subfigure}[b]{0.48\textwidth}
        \centering   \includegraphics[width=\textwidth]{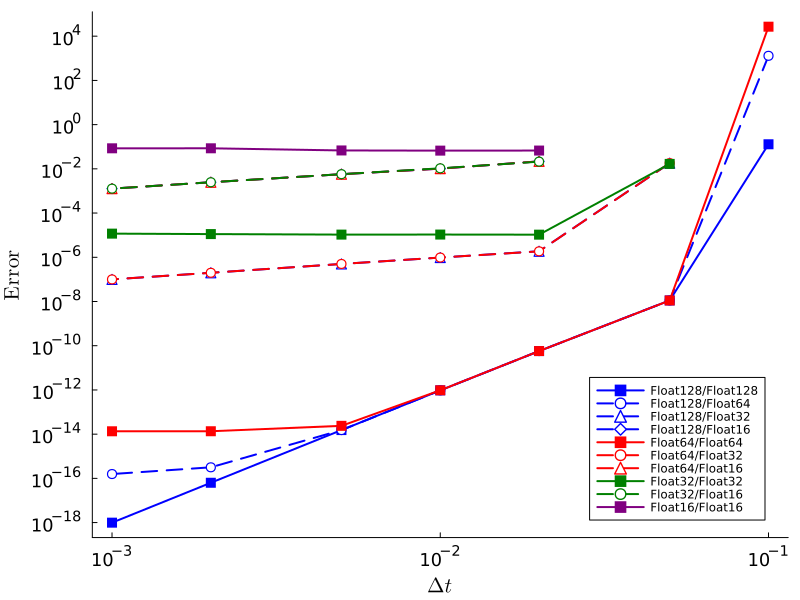}
    \end{subfigure}
    \caption{The mixed precision   sixth order  TDRK4s6p method \eqref{MP4s6p1e} 
    applied to Burgers' equation for (left) $N_x = 50$, and (right) $N_x = 100$.}\label{fig:burgers6p}
\end{figure}

Figures \ref{fig:burgers5p} and \ref{fig:burgers6p}
show the errors from  the mixed precision fifth order 
\eqref{Tsai3s5p1e} and sixth order \eqref{MP4s6p1e} methods 
(respectively) on the inviscid Burgers' equation above with  
$N_x =50$ points (left) and $N_x =100$ points (right).
These figures show mixed precision simulations with precisions from quad precision (\texttt{Float128}) to half precision (\texttt{Float16}).
For both the fifth order method in 
Figure \ref{fig:burgers5p}
and the sixth order method in 
Figure \ref{fig:burgers6p} we see that larger  $N_x$ and lower precision
result in a smaller $\dt$ required for stability. 
For fifth order, the $N_x =50$ plots (left) are stable
for double precision, quad precision, and mixed quad/double computation. 
For the mixed double/single and quad/single, the 
errors are high, for the largest $\dt$ tested, 
which suggests incipient instability.
For the half precision as well
as the mixed double/half and single/half the methods are clearly not stable for the largest $\dt$ tested.
The case for $N_x =100$ is clearer, as the methods that are below double precision are all unstable for the 
largest two values of $\dt$.
In all cases, the mixed quad/double codes show fifth order convergence until $\dt$ is small enough, at which point the perturbation error of order one becomes visible.
The value of $\dt$ for which this happens is smaller for 
$N_x=50$ than for $N_x =100$, which is expected since
the value of $\ep$ will depend not only on the precision but on the size of the matrix $\mD_x$ as well.

For sixth order, the $N_x =50$ plots (left) are all stable, while the 
$N_x =100$ plots (right)  are stable
for all $\dt$ for the double precision, quad precision, and mixed quad/double computation. A smaller 
$\dt $ is needed for stability for 
single precision and mixed quad/single and double/single. 
An even smaller $\dt $ is needed for stability for half precision, mixed  single/half, and double/half.
As in the fifth order plots in Figure 
\ref{fig:burgers5p}, we see that the mixed quad/double codes show sixth order convergence until $\dt$ is small enough, at which point the perturbation error of order one becomes visible. As expected,
the value of $\dt$ for which this happens is smaller for 
$N_x=50$ than for $N_x =100$.

\section{Conclusions} \label{sec:conclusions}
In this work we extended  Grant's additive framework for mixed accuracy and mixed precision Runge--Kutta methods 
\cite{Grant2022} to  mixed precision two-derivative methods. 
This enables us to naturally adopt a mixed precision computation 
of the second derivative into existing $p=3,4,5,6$th order methods \eqref{MP2s3p1e},
\eqref{MP2s4p1e}, \eqref{Tsai3s5p1e}, and \eqref{MP4s6p1e} (respectively),  and show using the perturbed order conditions that we expect a global order of
$O(\dt^p) + O(\ep \dt)$.  This framework further allowed us to develop third order methods 
that have higher order perturbation errors, the TDRK2s3p2e method \eqref{MP2s3p2e} with global order of
$O(\dt^3) + O(\ep \dt^2)$, and the  TDRK3s3p3e method \eqref{MP3s3p3e} with global order of
$O(\dt^3) + O(\ep \dt^3)$, and a fourth order TDRK3s4p2e method  \eqref{MP3s4p2e} with global order of
$O(\dt^4) + O(\ep \dt^2)$.  The numerical simulations in Section \ref{sec:NumResults}
confirm that the methods perform as expected in the time-evolution of  linear and nonlinear 
PDEs. 

The mixed accuracy framework we presented in Section \ref{sec:framework} can be 
generalized to mixed precision computations of $F$ as well as $\dot{F}$. It  can be
extended to higher order and adapted to perturbations that are smooth, or to higher orders in $\ep$
as well as $\dt$. 
Wherever two-derivative Runge--Kutta methods are used, this approach 
can be modified and extended to allow for significant computational savings.
Future directions for this approach include these extensions, as well as a formulation
of this framework for general linear methods (GLMs). Further developments will be 
application dependent,  in order to optimize for efficiency.  
For example, depending on the trade-off between the 
computational cost of the full precision operators and the  lower precision operators, 
it may be worthwhile to produce higher order methods that will use 
different balance of full precision and reduced precision computations.
These methods may be optimized for other properties including linear stability 
regions or strong stability properties. Now that the perturbation order conditions 
are determined, mixed precision two-derivative Runge--Kutta methods can be tailored and 
optimized as needed for different problems, as determined by the computational cost of 
each component.

\section*{Author Contribution Statement}

\noindent{\bf Sigal Gottlieb} was responsible for conceptualization of this project.
SG was responsible for developing the new fourth order methods in Section \ref{sec:MPTDRK4p}, and was primarily responsible for choosing the problems and methods in Section \ref{sec:NumResults}, and in writing the results in that section.
SG wrote the first draft of this manuscript, and was responsible for the final edits.

\noindent{\bf Zachary J. Grant} was responsible for developing the perturbation  order conditions, and for developing the new third order methods in Section \ref{sec:MPTDRK3p}.
ZJG modified LeVeque's stability region codes to produce all the plots in Section \ref{sec:MPTDRKmethods}.
ZJG co-wrote Sections \ref{sec:back1} and \ref{sec:back2}, carefully proofread the entire paper, and made numerous
editorial suggestions to improve the presentation.

\noindent{\bf Cesar Herrera}
was responsible for all the codes, computations, and preparation of graphs and tables in Section
\ref{sec:NumResults}, was involved in discussions on the presentation of the material. 
CH provided the expertise on mixed precision implementation in julia language, 
and all the numerical results. He read, commented, and edited the entire manuscript.

\bibliographystyle{plain}

\end{document}